\documentclass[a4paper,10pt]{article}

\usepackage{setspace}
\usepackage[bottom]{footmisc}
\usepackage{color}
\usepackage{subfig}
\usepackage{hyperref}
\usepackage{multirow}
\usepackage{cite}
\usepackage{amsmath} 
\usepackage{cleveref}
\usepackage{mathtools}
\usepackage{graphicx}
\usepackage{epsfig}
\usepackage{epstopdf}
\usepackage{float}
\usepackage{titling}
\usepackage[margin=1.0in]{geometry}
\usepackage{etoolbox}
 
\title{A compact subcell WENO limiting strategy using immediate neighbors for Runge-Kutta discontinuous Galerkin methods for unstructured meshes}
\makeatletter
\renewcommand\@date{{%
  \vspace{-\baselineskip}%
  \large\centering
  \begin{tabular}{@{}c@{}}
    S R Siva Prasad Kochi\textsuperscript{1} \\
    \normalsize siva.ksr@gmail.com
  \end{tabular}%
  \quad and\quad
  \begin{tabular}{@{}c@{}}
    M Ramakrishna\textsuperscript{2} \\
    \normalsize krishna@ae.iitm.ac.in
  \end{tabular}

  \bigskip

  \textsuperscript{1}Doctoral Candidate, Dept. of Aerospace Engg., IIT Madras.\par
  \textsuperscript{2}Professor, Dept. of Aerospace Engg., IIT Madras.

  \bigskip

  \today
}}
\makeatother
\begin{document}

\maketitle

\begin{abstract}
 In this paper, we generalize the compact subcell weighted essentially non oscillatory (CSWENO) limiting strategy for Runge-Kutta discontinuous Galerkin method developed recently in \cite{srspkmr1} for structured meshes to unstructured triangular meshes. The main idea of the limiting strategy is to divide the immediate neighbors of a given cell into the required stencil and to use a WENO reconstruction for limiting. This strategy can be applied for any type of WENO reconstruction. We used the WENO reconstruction proposed in \cite{zs4} and provided accuracy tests and results for two-dimensional Burgers' equation and two dimensional Euler equations to illustrate the performance of this limiting strategy.
 
\end{abstract}

\section{Introduction}\label{sec:intro}

In this paper, we look at the solution of hyperbolic conservation laws with the Runge-Kutta discontinuous Galerkin (RKDG) method \cite{cs6} on unstructured meshes. To control spurious oscillations near discontinuities, a limiter is used. For higher orders, weighted essentially non oscillatory (WENO) limiters are preferred as they maintain the order of the scheme. We generalize the WENO limiting strategy developed recently in \cite{srspkmr1} for structured meshes to unstructured triangular meshes. This strategy works with any type of WENO reconstruction in the target cell using only the immediate neighbors by dividing them into subcells to obtain the required stencil. We used the WENO reconstruction proposed in \cite{zs4} as it is quite simple for implementation.
\\
\\
\noindent WENO limiting for RKDG method was first presented by Qiu and Shu for structured meshes in \cite{qs1}. Their limiter was extended to unstructured meshes by Zhu et al. in \cite{zqsd}. This limiter requires neighbors of neighbors for limiting in a given cell for higher orders. Various other limiters given in \cite{zzsq1}, \cite{zzsq2}, \cite{zqz1} addressed this issue of using a wide stencil and used the polynomial in a given cell along with the polynomials of the immediate neighbors to obtain a limited polynomial using WENO reconstruction. On the other hand, Dumbser et al \cite{dl1} used a different strategy where the target cell is divided into subcells and an \textit{a posteriori} limiting strategy is used. This was further refined using an Adaptive Mesh Refinement (AMR) technique by Giri and Qiu in \cite{gq1}.
\\
\\
\noindent We extend the compact limiting strategy used in \cite{srspkmr1} to unstructured meshes. Here, we use the immediate neighbors for limiting in a given cell and these neighbors are divided into subcells for higher order limiting. We assign appropriate values to these subcells as explained in 
Section 3. This limiting strategy can be used with any type of WENO reconstructions as given in \cite{dk1} (called type-I WENO reconstruction), or \cite{zqsd} (called type-II WENO reconstruction) or \cite{lz1} (mixed reconstruction) or the more recent methods given in \cite{zq1} and \cite{zs4}. WENO schemes form a weighted combination of several local reconstructions based on different stencils (called small stencils) and use it as the final WENO reconstruction. Type-I reconstruction consists of WENO schemes whose order of accuracy is not higher than that of the reconstruction on each small stencil. These schemes require very wide stencils for higher orders. On the other hand, type-II WENO schemes require smaller stencils whose order of accuracy is higher than that of the reconstruction on each small stencil. Each of these reconstructions have their pros and cons and they are discussed in detail in \cite{lz1}. We use the WENO reconstruction given by Zhu and Shu in \cite{zs4} as it is quite simple in implementation and extension to higher orders is easy. These \cite{zs4} schemes have linear weights that can be any positive numbers on the condition that their sum is one and number of spatial stencils is smaller than that of the same order accurate classical finite volume WENO schemes \cite{hs1} on triangular meshes. We emphasize that our procedure works with any of the WENO reconstructions listed above by dividing the neighbors of the cell in which limiting is needed as is required by the reconstruction. We call this limiting strategy as compact subcell WENO limiter or CSWENO limiter in short.
\\
\\
\noindent The paper is organized as follows. We describe the formulation of the discontinuous Galerkin method used for all our results in 
Section 2, the proposed limiting procedure is described in 
Section 3 and the testing of the limiter and the results are described in 
Section 4 and finally we conclude the paper in 
Section 5.

\section{Formulation of discontinuous Galerkin method on unstructured meshes}\label{sec:formulation}

\noindent Consider a two-dimensional conservation law of the form for variable $u$

\begin{equation}\label{governEqn2D}
 \frac{\partial u}{\partial t} + \frac{\partial f(u)}{\partial x} + \frac{\partial g(u)}{\partial y} = 0, \qquad \qquad x,y\in \mathbf{D}
\end{equation}

\noindent with the initial condition at $t=0$ given by,

\begin{displaymath}
 u(x,y,0) = u_{0}(x,y),
\end{displaymath}

\noindent and the appropriate boundary conditions on the boundary $\partial \mathbf{D}$. Here $f(u)$ and $g(u)$ are fluxes in $x$ and $y$ directions respectively.
\\
\\
\noindent We assume that the domain $\mathbf{D}$ can be triangulated by $K$ elements as

\begin{equation}\label{triangulation}
 \mathbf{D} = \bigcup_{k=1}^{K} \mathbf{I}_{k}
\end{equation}

\noindent where $\mathbf{I}_{k}$ is a straight sided triangle and the triangulation is assumed to be geometrically conforming i.e., $\partial \mathbf{D}$ is approximated by a piece wise linear polygon with each line a face of the triangle.
\\
\\
\noindent We approximate the local solution as a polynomial of order $N$ as given by:

\begin{equation}\label{modalForm2D}
 u_{h}^{k}(x,y,t) = \sum_{n=0}^{N_{p}-1} \hat{u}^{k}_{n}(t)\psi_{n}^{k}(x,y) \qquad \forall x,y\in\mathbf{I}_{k}
\end{equation}

\noindent where $\psi_{n}^{k}(x,y)$ is a two dimensional polynomial basis of order $N$, $N_{p}$ is the number of degrees of freedom and $h$ is the size of the grid (which is in general the average of the lengths of the sides of the triangle). This is termed to be $\mathbf{P}^{N}$ based discontinuous Galerkin method whose formal order of accuracy is $N+1$. The number of degrees of freedom $N_{p}$ is given by:

\begin{displaymath}
 N_{p} = \frac{(N+1)(N+2)}{2}
\end{displaymath}
\\
\\
\noindent We use the orthonormal basis as given in \cite{hestha1} for a standard isoparametric triangle $I$ such that

\begin{equation}\label{triangleDomain}
 I = \{\vec{r} = (r,s)|(r,s)\geq -1 ; r+s \leq 0\}
\end{equation}
\\
\\
\noindent The basis is given as:

\begin{equation}\label{triangleBasis}
 \psi_{m}(\vec{r}) = \sqrt{2} \left[P_{i}(a)P_{j}^{(2i+1,0)}(b)\right](1-b)^{i}
\end{equation}

\noindent where

\begin{displaymath}
 m = j+(N+1)i+1-\frac{i}{2}(i-1),\qquad (i,j)\geq 0; \qquad i+j\leq N
\end{displaymath}

\noindent and 

\begin{displaymath}
 a = 2\frac{1+r}{1-s} - 1;\qquad b = s
\end{displaymath}

\noindent and $P_{n}^{(\alpha,\beta)}(x)$ is the $n$th order Jacobi polynomial.
\\
\\
\noindent We approximate the fluxes $f(u(x,y,t),x,y,t)$ and $g(u(x,y,t),x,y,t)$ in the domain $\mathbf{D}$ as

\begin{equation}\label{fluxApprox1}
 f_{h}^{k}(u_{h}^{k}) = \sum_{n=0}^{N_{p}-1} \hat{f}^{k}_{n}(t)\psi_{n}^{k}(x,y) \qquad \forall x,y\in\mathbf{I}_{k}
\end{equation}

\begin{equation}\label{fluxApprox2}
 g_{h}^{k}(u_{h}^{k}) = \sum_{n=0}^{N_{p}-1} \hat{g}^{k}_{n}(t)\psi_{n}^{k}(x,y) \qquad \forall x,y\in\mathbf{I}_{k}
\end{equation}

\noindent where $\hat{f}^{k}_{n}(t) = f(\hat{u}^{k}_{n}(t))$ and $\hat{g}^{k}_{n}(t) = g(\hat{u}^{k}_{n}(t))$. Substituting equations \eqref{modalForm2D}, \eqref{fluxApprox1} and \eqref{fluxApprox2} in \eqref{governEqn2D} and integrating it by parts, we get the following scheme to advance the degrees of freedom $\hat{u}^{k}_{n}(t)$ in time in an element given by $\mathbf{I}_{k}=[x_{l}^{k},y_{l}^{k}]\times [x_{r}^{k},y_{r}^{k}]$:

\begin{align}\label{weakFormScheme2D}
 \frac{d}{dt}\hat{u}_{h}^{k} = (\mathbf{M}^{k})^{-1}(\mathbf{S_{x}}^{k})^{T} \hat{f}_{h}^{k}(\hat{u}_{h}^{k}) + (\mathbf{M}^{k})^{-1}(\mathbf{S_{y}}^{k})^{T} \hat{g}_{h}^{k}(\hat{u}_{h}^{k}) \\ \nonumber - (\mathbf{M}^{k})^{-1} (f^{*}|_{r_{N_{p}}}e_{N_{p}} - f^{*}|_{r_{1}}e_{1}) - (\mathbf{M}^{k})^{-1} (g^{*}|_{s_{N_{p}}}e_{N_{p}} - g^{*}|_{s_{1}}e_{1}) = L(u_{h}^{k})
\end{align}

\noindent Here, $\hat{u}_{h}^{k} = [\hat{u}_{0}^{k} \ldots \hat{u}_{N_{p}-1}^{k}]^{T}$, $e_{i}$ is a vector of dimension $N_{p}$ which has zero entries everywhere except at the $i$th location and $\mathbf{M}^{k}$ is the local mass matrix which is given as:


\begin{equation}\label{massMatrix2D}
 \mathbf{M}^{k} = \left[M_{ij}^{k}\right] = \left[\int_{\mathbf{I}_{k}} \psi_{i}^{k}(x,y) \psi_{j}^{k}(x,y) ~\text{dA}\right]
\end{equation}

\noindent and $\mathbf{S_{x}}^{k}$ and $\mathbf{S_{y}}^{k}$ are the local stiffness matrices given by:

%

\begin{equation}\label{stiffnessMatrix12D}
 \mathbf{S_{x}}^{k} = \left[S_{xij}^{k}\right] = \left[\int_{\mathbf{I}_{k}} \psi_{i}^{k}(x,y) \frac{\partial \psi_{j}^{k}(x,y)}{\partial x} ~\text{dA}\right]
\end{equation}

\begin{equation}\label{stiffnessMatrix22D}
 \mathbf{S_{y}}^{k} = \left[S_{yij}^{k}\right] = \left[\int_{\mathbf{I}_{k}} \psi_{i}^{k}(x,y) \frac{\partial \psi_{j}^{k}(x,y)}{\partial y} ~\text{dA}\right]
\end{equation}

\noindent Also, $f^{*}$ and $g^{*}$ are the $x$ and $y$ components of the monotone numerical flux at the interface which is calculated using an exact or approximate Riemann solver. We have used the Lax-Friedrichs flux for all the test cases given below.
\\
\\
\section{Proposed limiting procedure on unstructured meshes}\label{sec:limiter}

\noindent In this section, we describe the details of the proposed limiting procedure using WENO reconstruction and the division of cells into subcells for the discontinuous Galerkin method. This step is a generalization of the procedure in \cite{srspkmr1} for structured meshes. The common method for limiting in discontinuous Galerkin method is:
\\
\noindent \textbf{1)} Identify the cells which need to be limited, known as troubled cells. \\
\noindent \textbf{2)} Replace the solution polynomial in the troubled cell with a new polynomial that is less oscillatory but with the same cell average and order of accuracy.
\\
\\
\noindent For the first step, we have used the KXRCF troubled cell indicator \cite{kxrcf} for all the calculations done in this paper as it is rated highly by Qiu and Shu in \cite{qs2} on the basis of it's performance in detecting the discontinuities in various test problems. We give a brief description of the troubled cell indicator in the following subsection.

\subsection{KXRCF Troubled cell indicator}

\noindent The troubled cell indicator used in this work was developed by Krivodonova et al. \cite{kxrcf}. This indicator is termed as KXRCF troubled cell indicator by using the names of the authors of the paper \cite{kxrcf}.
\\
\\
\noindent For a given problem, partition the boundary $\partial\mathbf{I}_{j}$ of a given cell $\mathbf{I}_{j}$ into portions $\partial\mathbf{I}_{j}^{-}$ and $\partial\mathbf{I}_{j}^{+}$ where the flow is into ($\bar{v}.\bar{n}<0$) and out of ($\bar{v}.\bar{n}>0$) $\mathbf{I}_{j}$, respectively. For the scalar conservation laws, $\bar{v}$ is taken to be $(f'(u),g'(u))$ and for the Euler equations it is the velocity vector. The troubled cell indicator is defined as

\begin{equation}\label{troubledCellIndicator}
R_{j} = \frac{|\int_{\partial\mathbf{I}_{j}^{-}} (U_{j}-U_{nbj}) ds|}{h^{(N+1)/2}|\partial\mathbf{I}_{j}^{-}|\quad ||U_{j}||}
\end{equation}

\noindent where $U_{j}$ is the discontinuous Galerkin value of $u$ on $\mathbf{I}_{j}$, $U_{nbj}$ the value across the boundary, $||U_{j}||$ the standard $L^{2}$ norm in the cell $\mathbf{I}_{j}$  and $N$ is the order of the polynomial basis. We take $h$ to be the radius of the circumscribed circle in $\mathbf{I}_{j}$. Now if $R_{j}>C_{K}$, $\mathbf{I}_{j}$ is identified to be a troubled cell. We take $C_{K}$ to be 1.
\\
\\
\subsection{Limiting procedure}

\begin{figure}[H]
\begin{center}
\includegraphics[scale=0.5]{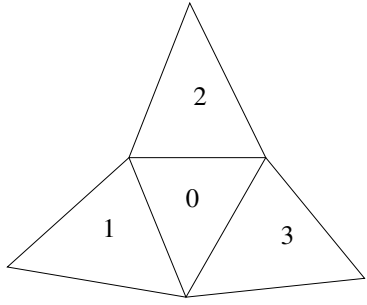}
\caption{Troubled Cell and its immediate neighbors}
\label{fig1:originalStencil}
\end{center}
\end{figure}

\noindent After identifying the troubled-cells, we would like to reconstruct the values of $\hat{u}_{n}^{j}$ for the troubled-cell $\mathbf{I}_{j}$ for $n=1,\ldots,N_{p}-1$. That is, we retain the cell average $\hat{u}_{0}^{j}$ and reconstruct all the other degrees of freedom. To do that, we use the quadrature points given in \cite{wsj} (SCP-quadrature) and find the solution at those points using WENO reconstruction. These quadrature points are termed as $(x_{q},y_{q})$ for $q=0,\ldots,N_{p}-1$. We can use any of the WENO reconstructions given in \cite{dk1} (called type I WENO reconstruction), or \cite{zqsd} (called type II WENO reconstruction) or \cite{lz1} (mixed reconstruction) or the more recent methods given in \cite{zq1} and \cite{zs4}. WENO schemes form a weighted combination of several local reconstructions based on different stencils (called small stencils) and use it as the final WENO reconstruction. Type-I reconstruction consists of WENO schemes whose order of accuracy is not higher than that of the reconstruction on each small stencil. These schemes require very wide stencils for higher orders. On the other hand, type-II WENO schemes require smaller stencils whose order of accuracy is higher than that of the reconstruction on each small stencil. But, they sometimes require negative linear weights and do not work properly for geometries of poor mesh quality. Each of these reconstructions have their pros and cons and they are discussed in detail in \cite{lz1}. We use the WENO reconstruction given by Zhu and Shu in \cite{zs4} as it is quite simple in implementation and extension to higher orders is easy. These \cite{zs4} schemes have linear weights that can be any positive numbers on the condition that their sum is one and number of spatial stencils is smaller than that of the same order accurate classical finite volume WENO schemes \cite{hs1} on triangular meshes. We emphasize that our procedure works with any of the WENO reconstructions listed above by dividing the neighbors of the cell in which limiting is needed as is required by the reconstruction. For the $\mathbf{P}^{1}$ based DGM, we now describe the procedure for the reconstruction of the moments $\hat{u}_{1}^{k}$ and $\hat{u}_{2}^{k}$ in the troubled cell $\mathbf{I}_{j}$ using the WENO reconstruction procedure for triangles. We relabel the troubled cell as $\Delta_{0}$ and its immediate neighbors as $\Delta_{1},\Delta_{2}$ and $\Delta_{3}$ as shown in Figure \ref{fig1:originalStencil} where they are labeled according to their subscripts.

\begin{figure}[htbp]
  \centering
  \subfloat[Division of cell neighbors for $\mathbf{P}^{1}$ based DGM]{\label{fig:NeighborP1}\includegraphics[width=0.45\textwidth]{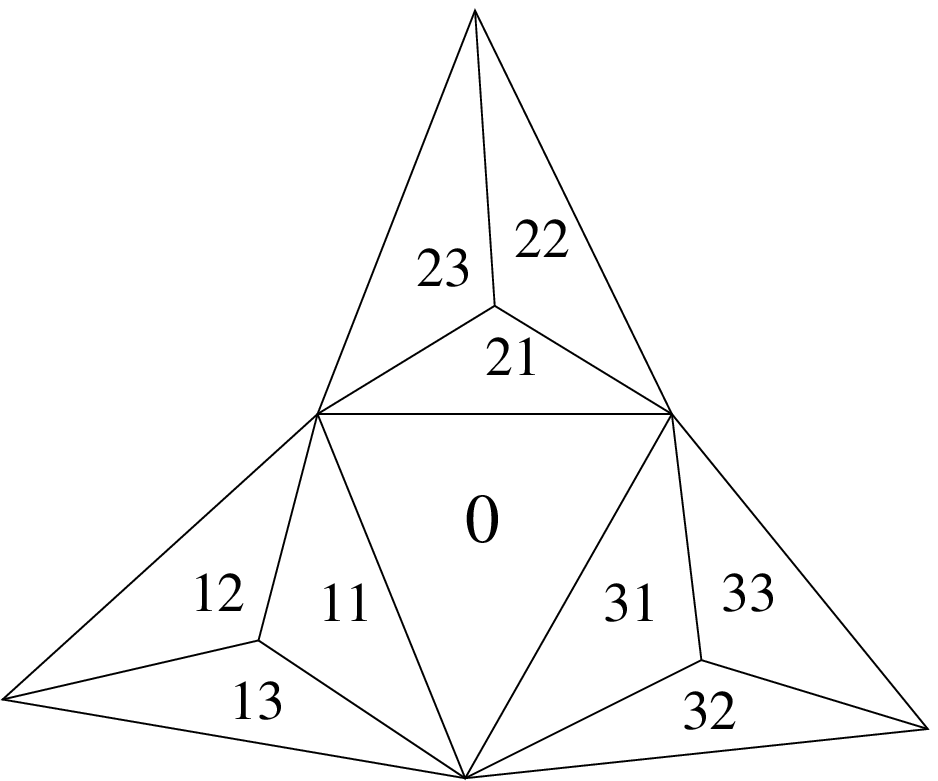}}
  \subfloat[Division of cell neighbors for $\mathbf{P}^{2}$ based DGM]{\label{fig:NeighborP2}\includegraphics[width=0.45\textwidth]{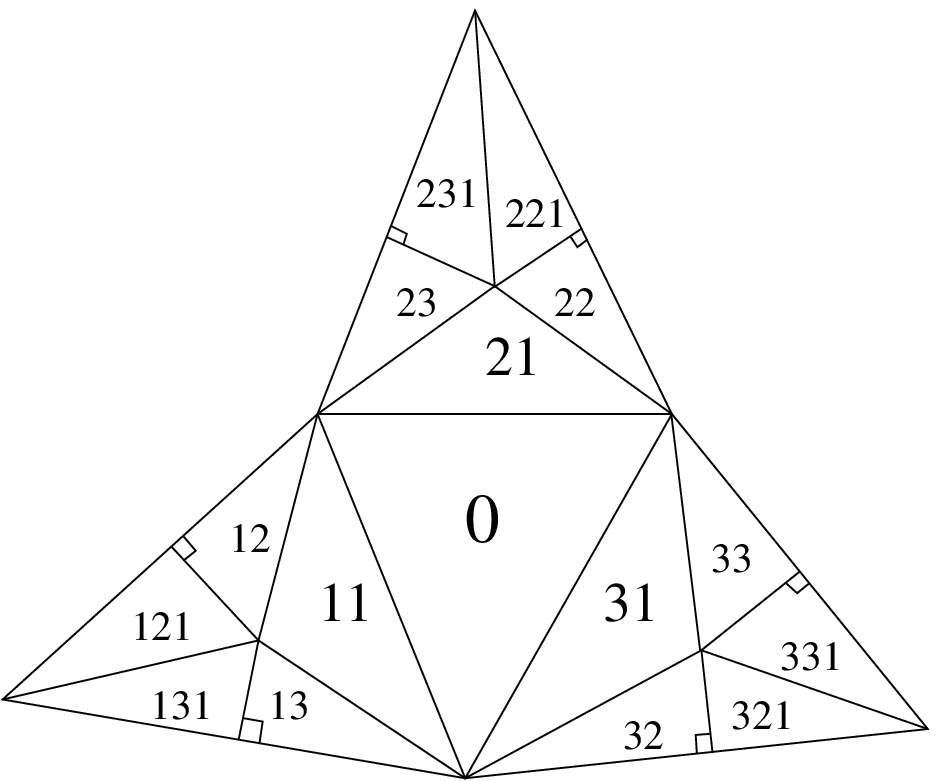}}
  \caption{Division of cell neighbors for WENO reconstruction using $\mathbf{P}^{1}$ and $\mathbf{P}^{2}$ based DGM}
  \label{fig:NeighborDivision}
\end{figure}

\noindent For a third order WENO reconstruction, we also need the neighbors of the cells $\Delta_{1},\Delta_{2}$ and $\Delta_{3}$ as suggested in \cite{zqsd}. Instead of using their actual neighbors, we divide each of those triangles into three equal parts by joining the centroid of a triangle to its vertices as shown in Figure \ref{fig:NeighborP1}. Now, we get the new stencil required for the third order WENO reconstruction. Now, we use the DG local polynomial for $\Delta_{1}$ and take its average in the new cells $\Delta_{11}$, $\Delta_{12}$ and $\Delta_{13}$ to obtain the cell averages for the new cells. Similarly, we can calculate the cell average for the cells $\Delta_{21}$, $\Delta_{22}$, $\Delta_{23}$, $\Delta_{31}$, $\Delta_{32}$ and $\Delta_{33}$ using the DG local polynomial of cells $\Delta_{2}$ and $\Delta_{3}$. We can use the same procedure to divide the neighboring cells for $\mathbf{P}^{N}$ based DGM for any $N$ based on the required number of quadrature points. As the order of the polynomial increases, the number of triangles required to maintain the order of the scheme increases. We split the neighbors in such a way that this will not reduce the order of accuracy as explained in \cite{hs1}.

\begin{figure}[H]
\begin{center}
\includegraphics[scale=0.5]{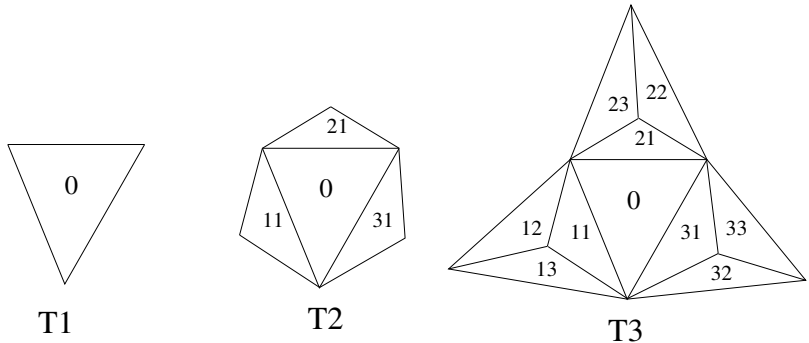}
\caption{Required stencils for $\mathbf{P}^{1}$ based DGM}
\label{fig:StencilsP1}
\end{center}
\end{figure}

\noindent Now the step by step WENO reconstruction procedure to be used for $\mathbf{P}^{1}$ based DGM is described below:
\\
\\
\noindent \textbf{Step 1:} We use stencils T1, T2 and T3 as shown in Figure \ref{fig:StencilsP1} for WENO reconstruction. For T1 = \{$\Delta_{0}$\}, we assign a zeroth order polynomial $q_{1}(x,y)=\hat{u}_{0}^{j}$. Now on T2 = \{$\Delta_{0}$,$\Delta_{11}$,$\Delta_{21}$,$\Delta_{31}$\}, we obtain a linear polynomial $q_{2}(x,y)$ such that it has the same cell average as $u$ on the troubled cell $\Delta_{0}$ and matches the cell averages of $u$ (obtained earlier) on the other triangular cells in the set T2$\backslash$\{$\Delta_{0}$\} in a least square sense \cite{hs1}. To obtain a quadratic polynomial $q_{3}(x,y)$, we use the stencil T3 = \{$\Delta_{0}$,$\Delta_{11}$,$\Delta_{21}$,$\Delta_{31}$,$\Delta_{12}$,$\Delta_{13}$,$\Delta_{22}$,$\Delta_{23}$,$\Delta_{32}$,$\Delta_{33}$\} such that it has the same cell average as $u$ on the troubled cell $\Delta_{0}$ and matches the cell averages of $u$ on the other triangular cells in the set T3$\backslash$\{$\Delta_{0}$\} in a least square sense \cite{hs1}.
\\
\\
\noindent \textbf{Step 2:} Now, we obtain equivalent expressions for these polynomials $q_{1}(x,y)$, $q_{2}(x,y)$ and $q_{3}(x,y)$. We define:
\begin{equation}\label{equivalentOrder0}
 p_{1}(x,y)=q_{1}(x,y)
\end{equation}

\begin{equation}\label{equivalentOrder1}
 p_{2}(x,y) = \frac{1}{\gamma_{2,2}}q_{2}(x,y) - \frac{\gamma_{1,2}}{\gamma_{2,2}}p_{1}(x,y) \quad \text{with} \quad \sum_{l=1}^{2}\gamma_{l,2}=1 \quad \text{and} \quad \gamma_{2,2} \neq 0
\end{equation}

\noindent and

\begin{equation}\label{equivalentOrder2}
 p_{3}(x,y) = \frac{1}{\gamma_{3,3}}q_{3}(x,y) - \sum_{l=1}^{2}\frac{\gamma_{l,3}}{\gamma_{3,3}}p_{l}(x,y) \quad \text{with} \quad \sum_{l=1}^{3}\gamma_{l,3}=1 \quad \text{and} \quad \gamma_{3,3} \neq 0
\end{equation}

\noindent Here the $\gamma_{l,m}$ are linear weights. We choose them as in \cite{zs4} as $\gamma_{l,m} = \bar{\gamma}_{l,m}/\left({\sum_{l=1}^{m}\bar{\gamma}_{l,m}}\right)$ with $\bar{\gamma}_{l,m}=(1/h)^{l-1}$, $l=1,\ldots , t$ and $m=2,3$ ($t=m$) where $h$ is the grid size (which is in general the average of the lengths of the sides of the triangle). Here we changed the expression for $\bar{\gamma}_{l,m}$ from $\bar{\gamma}_{l,m}=10^{l-1}$ (which is used in \cite{zs4}) to $\bar{\gamma}_{l,m}=(1/h)^{l-1}$. This change was made so that $\gamma_{l,m}$ = $O(h^{m-l})$. In this manner, we get larger linear weights for high degree polynomials and smaller linear weights for low degree polynomials. Also, using a Taylor series expansion of cell averages about the barycenter of the triangle $\Delta_{0}$ as given in \cite{sonar1}, we obtain $u(x_{q},y_{q})-p_{l}(x_{q},y_{q})=O(h^{l})$, where $u(x_{q},y_{q})$ is the exact solution at the quadrature point $q$.
\\
\\
\noindent \textbf{Step 3:} We now compute the smoothness indicator for each stencil denoted by $\beta_{l}$ for $l = 2,3$. As given by \cite{js1}, we use:

\begin{equation}\label{smoothnessIndicator}
\beta_{l} = \sum_{|l|=1}^{N}|\Delta_{0}|^{|l|-1}\int_{\Delta_{0}} \left(\frac{\partial^{|l|}}{\partial x^{l_{1}}\partial x^{l_{2}}}p_{i}(x,y)\right)^{2} dxdy, \quad l = 2,3
\end{equation}

\noindent where $l=(l_{1},l_{2})$, $|l|=l1+l2$. For $\beta_{1}$, we magnify it from 0 to a value defined in equation \eqref{smoothnessIndicatorT1} as given by \cite{zs4}. Using the stencils in Figure \ref{fig:StencilsP1}, we construct three polynomials $p_{1,l}(x,y) \in \text{span}\{\frac{x-x_{l}}{|\Delta_{l}|^{1/2}}\}$, $l=1,2,3$, satisfying $p_{1,1}(x_{12},y_{12})=\bar{u}_{12}-\bar{u}_{11}$, $p_{1,1}(x_{13},y_{13})=\bar{u}_{13}-\bar{u}_{11}$, $p_{1,2}(x_{22},y_{22})=\bar{u}_{22}-\bar{u}_{21}$, $p_{1,2}(x_{23},y_{23})=\bar{u}_{23}-\bar{u}_{21}$, $p_{1,3}(x_{32},y_{32})=\bar{u}_{32}-\bar{u}_{31}$, $p_{1,3}(x_{33},y_{33})=\bar{u}_{33}-\bar{u}_{31}$, where $(x_{l},y_{l})$ are centroids of $\Delta_{l}$, $l = 11,12,13;21,22,23;31,32,33$, respectively. Now, we use equation \eqref{smoothnessIndicator} to find $\beta_{1,l}$, $l=1,2,3$. We set $\lambda_{1,1}=\lambda_{1,2}=\lambda_{1,3}=1/3$ and

\begin{equation}\label{sigmaEqn}
 \sigma_{l}=\lambda_{1,l}\left(1+\frac{\left(\left[|\beta_{1,1}-\beta_{1,2}|^{2}+|\beta_{1,2}-\beta_{1,3}|^{2}+|\beta_{1,3}-\beta_{1,1}|^{2}\right]/3\right)^{2}}{\beta_{1,l}+\epsilon}\right),l=1,2,3
\end{equation}

\noindent with $\sigma=\sigma_{1}+\sigma_{2}+\sigma_{3}$, where $\epsilon=10^{-10}$. Then we get

\begin{equation}\label{smoothnessIndicatorT1}
\beta_{1} = \sum_{|l|=1}|\Delta_{0}|^{|l|} \left(\frac{\partial^{|l|}}{\partial x^{l_{1}}\partial x^{l_{2}}}\left(\sum_{l=1}^{3} \frac{\sigma_{l}}{\sigma} p_{1,l}(x,y)\right)\right)^{2}
\end{equation}

\noindent Again using a Taylor series expansion of cell averages about the barycenter of the triangle $\Delta_{0}$, we obtain 
for $l=1,2,3$, $\beta_{l} = d(x_{l},y_{l})h^{2}(1+O(h^{2}))$ for some function $d$ where $(x_{l},y_{l})$ is the location of barycenter of $\mathbf{I}_{j}$. If a particular stencil, say T2 is not smooth, then $\beta_{2}=O(1)$.
\\
\\
\noindent \textbf{Step 4:} We compute the nonlinear weights using the WENO-Z recipe \cite{ccd}. We first find $\tau$ using

\begin{equation}\label{tauNonlinearWeightsP1}
 \tau = \frac{|\beta_{3}-\beta_{1}|^{2}+|\beta_{3}-\beta_{2}|^{2}}{4}
\end{equation}

\noindent Here, using the Taylor series expansions about the barycenter of the triangle $\Delta_{0}$, we obtain $\tau=O(h^{6})$. Now, the nonlinear weights are given by

\begin{equation}\label{nonLinearWeightsP1}
\omega_{l}=\frac{\bar{\omega}_{l}}{\sum_{l=1}^{3}\bar{\omega}_{l}}, \quad \bar{\omega}_{l} = \gamma_{l,3}\left(1+\frac{\tau}{\epsilon + \beta_{l}}\right) \quad l = 1,2,3
\end{equation}

\noindent Again we take $\epsilon$ as $10^{-10}$. Now assuming $\epsilon << \beta_{l}$, we get $\omega_{1}-\gamma_{1,3}=O(h^{6})$, $\omega_{2}-\gamma_{2,3}=O(h^{5})$, and $\omega_{3}-\gamma_{3,3}=O(h^{4})$. If a particular stencil, say T2 is not smooth, then $\omega_{2}-\gamma_{l,2}=O(1)$. We note that the stencil T1 is always smooth as it contains only one cell with a given cell average. Now, we can write the final WENO approximation at the required quadrature point $q$ as

\begin{equation}\label{WENOApproxQuadratureP1}
u_{q} = \sum_{l=1}^{3} \omega_{l} p_{l}(x_{q},y_{q}) \quad q=0,\ldots,N_{p}-1
\end{equation}

\noindent Now, we write

\begin{equation}\label{orderApprox}
 u_{q} - u(x_{q},y_{q}) = \omega_{0}p_{0}(x_{q},y_{q}) + \omega_{1}p_{1}(x_{q},y_{q}) + \omega_{2}p_{2}(x_{q},y_{q}) - u(x_{q},y_{q})
\end{equation}

\noindent which can written as

\begin{align}\label{orderApproxExpanded}
 u_{q} - u(x_{q},y_{q}) = [(\omega_{3} - \gamma_{3,3})(p_{3}(x_{q},y_{q})-u(x_{q},y_{q}))] + [(\omega_{2} - \gamma_{2,3})(p_{2}(x_{q},y_{q})-u(x_{q},y_{q}))] \nonumber \\ 
  + [(\omega_{1} - \gamma_{1,3})(p_{1}(x_{q},y_{q})-u(x_{q},y_{q}))] + [\gamma_{3,3}(p_{3}(x_{q},y_{q})-u(x_{q},y_{q}))] \nonumber \\ + [\gamma_{2,3}(p_{2}(x_{q},y_{q})-u(x_{q},y_{q}))] + [\gamma_{1,3}(p_{3}(x_{q},y_{q})-u(x_{q},y_{q}))]
\end{align}

\noindent If all the stencils are smooth, from the right hand side of equation \eqref{orderApproxExpanded}, we can see that $u(x_{q},y_{q}) - u_{q} = O(h^{3})$. If any one of the stencils (say T2) is not smooth, then $\omega_{2}-\gamma_{2,3}=O(1)$ and the right hand side of equation \eqref{orderApproxExpanded} will become $O(h^{2})$. This limit is the required order of accuracy for the reconstruction of moments for $\mathbf{P^{1}}$ based DGM.

\noindent \textbf{Step 5:} Finally, we obtain the reconstructed degrees of freedom based on the reconstructed point values $u_{q}$ at the quadrature points $(x_{q},y_{q})$ and a numerical integration as

\begin{equation}\label{finalMoments}
\hat{u}_{i}^{j} = |\Delta_{0}| \sum_{q} w_{q} u_{q} \psi_{i}^{j}(x_{q},y_{q}) \quad i=1,\ldots,N_{p}-1
\end{equation}

\noindent where $|\Delta_{0}|$ is the area of $\Delta_{0}$ and $w_{G}$'s are the quadrature weights for the points $(x_{q},y_{q})$ as given by \cite{wsj}.
\\
\\
\noindent For the $\mathbf{P}^{2}$ based DGM, the procedure to construct the second order moments $\hat{u}_{1}^{k}$, $\hat{u}_{2}^{k}$, $\hat{u}_{3}^{k}$, $\hat{u}_{4}^{k}$ and $\hat{u}_{5}^{k}$ in the troubled cell $\Delta_{0}$ is similar to the above procedure. Here, we require some extra neighbors for the WENO reconstruction as we need a fifth order WENO reconstruction and they are shown in Figure \ref{fig:NeighborP2}. Here note that the new triangles $221,231$ etc., are obtained by dropping a perpendicular to the opposite side as required. The cell averages for the new cells are obtained using the same procedure described above for $\mathbf{P}^{1}$ based DGM. The WENO reconstruction procedure for this case is described below:
\\
\\
\begin{figure}[H]
\begin{center}
\includegraphics[scale=0.5]{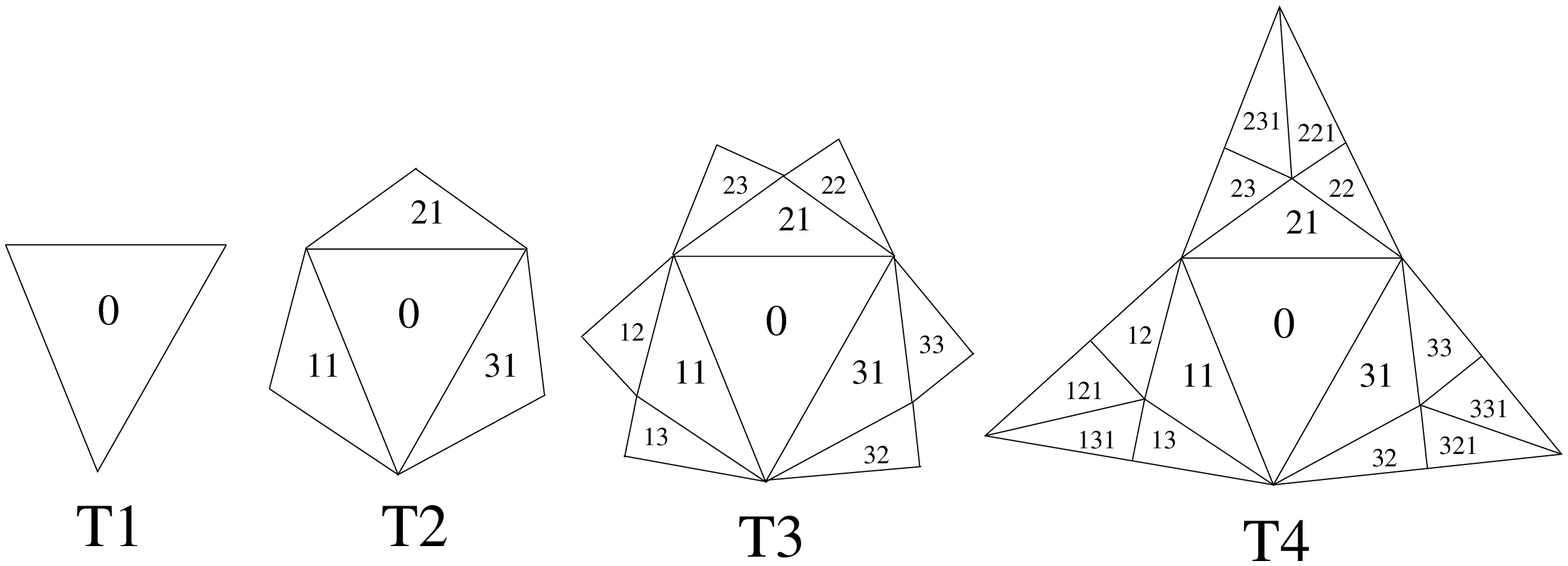}
\caption{Required stencils for $\mathbf{P}^{2}$ based DGM}
\label{fig:StencilsP2}
\end{center}
\end{figure}

\noindent \textbf{Step 1:} We use stencils T1, T2, T3 and T4 as shown in Figure \ref{fig:StencilsP2} for WENO reconstruction. We obtain $q_{1}(x,y)$, $q_{2}(x,y)$ and $q_{3}(x,y)$ using stencils T1, T2 and T3 as described in the procedure above. Now, using T4 = \{$\Delta_{0}$,$\Delta_{11}$,$\Delta_{21}$,$\Delta_{31}$,$\Delta_{12}$,$\Delta_{13}$,$\Delta_{121}$,$\Delta_{131}$,$\Delta_{22}$,$\Delta_{23}$,$\Delta_{221}$,$\Delta_{231}$,$\Delta_{32}$,$\Delta_{33}$,$\Delta_{321}$,$\Delta_{331}$\}, we obtain a cubic polynomial $q_{4}(x,y)$, such that it has the same cell average as $u$ on the troubled cell $\Delta_{0}$ and matches the cell averages of $u$ on the other triangular cells in the set T4$\backslash$\{$\Delta_{0}$\} in a least square sense. Again a quartic polynomial $q_{5}(x,y)$ is obtained on T5 such that it has the same cell average as $u$ on the troubled cell $\Delta_{0}$ and matches the cell averages of $u$ on the other triangular cells in the set T4$\backslash$\{$\Delta_{0}$\} in a least square sense.
\\
\\
\noindent \textbf{Step 2:} We obtain $p_{1}(x,y)$, $p_{2}(x,y)$ and $p_{3}(x,y)$ using equations \eqref{equivalentOrder0}, \eqref{equivalentOrder1} and \eqref{equivalentOrder2} respectively. We now define:

\begin{equation}\label{equivalentOrder3}
 p_{4}(x,y) = \frac{1}{\gamma_{3,4}}q_{4}(x,y) - \sum_{l=1}^{2}\frac{\gamma_{l,4}}{\gamma_{3,4}}p_{l}(x,y) \quad \text{with} \quad \sum_{l=1}^{3}\gamma_{l,4}=1 \quad \text{and} \quad \gamma_{3,4} \neq 0
\end{equation}

\begin{equation}\label{equivalentOrder4}
 p_{5}(x,y) = \frac{1}{\gamma_{4,5}}q_{5}(x,y) - \sum_{l=1}^{3}\frac{\gamma_{l,5}}{\gamma_{4,5}}p_{l}(x,y) \quad \text{with} \quad \sum_{l=1}^{4}\gamma_{l,5}=1 \quad \text{and} \quad \gamma_{4,5} \neq 0
\end{equation}

\noindent We choose $\gamma_{l,m}$ again as given above  for $m=2,3$. For $m=4,5$, we set $\gamma_{l,m} = \bar{\gamma}_{l,m}/\left(\sum_{l=1}^{m}\bar{\gamma}_{l,m}\right)$ with $\bar{\gamma}_{l,m}=(1/h)^{l-1}$, $l=1,\ldots , t$ and $m=4,5$ ($t=m-1$) where $h$ is the grid size.
\\
\\
\noindent \textbf{Step 3:} Using the same procedure given in Step 3 for $\mathbf{P}^{1}$ based DGM and equations \eqref{smoothnessIndicator}, \eqref{smoothnessIndicatorT1}, we obtain $\beta_{l}$, for $l=1,2,3,4,5$.
\\
\\
\noindent \textbf{Step 4:} We find $\tau$ required for calculating nonlinear weights using

\begin{equation}\label{tauNonlinearWeightsP2}
 \tau = \frac{|\beta_{4}-\beta_{1}|^{2}+|\beta_{4}-\beta_{2}|^{2}+|\beta_{4}-\beta_{3}|^{2}}{16}
\end{equation}

\noindent Now, we calculate the nonlinear weights using 

\begin{equation}\label{nonLinearWeightsP2}
\omega_{l}=\frac{\bar{\omega}_{l}}{\sum_{l=1}^{4}\bar{\omega}_{l}}, \quad \bar{\omega}_{l} = \gamma_{l,5}\left(1+\frac{\tau}{\epsilon + \beta_{l}}\right) \quad l = 1,2,3,4
\end{equation}

\noindent Again we take $\epsilon$ as $10^{-10}$. Now, we can write the final WENO approximation at the required quadrature point $q$ as

\begin{equation}\label{WENOApproxQuadratureP2}
u_{q} = \sum_{l=1}^{4} \omega_{l} p_{l}(x_{q},y_{q}) \quad q=0,\ldots,N_{p}-1
\end{equation}

\noindent Using equations \eqref{WENOApproxQuadratureP2} and \eqref{orderApproxExpanded}, we can say that when the solution of smooth, $u(x_{q},y_{q}) - u_{q} = O(h^{5})$. If one of the stencils is not smooth, we have $u(x_{q},y_{q}) - u_{q} = O(h^{3})$     . This limit is the required order of accuracy for the reconstruction of moments for $\mathbf{P^{2}}$ based DGM. We can reconstruct the degrees of freedom in the same way by dividing the neighbors and using WENO reconstruction for $\mathbf{P^{3}}$ based DGM.
\\
\\
Now, equation \eqref{finalMoments} gives us the degrees of freedom for $\mathbf{P}^{2}$ based DGM. For a non-orthonormal basis, we define 
\begin{equation}
D_{i}^{j} = |\Delta_{0}| \sum_{q} w_{q} u_{q} \psi_{i}^{j}((x_{q},y_{q})) \quad i=1,\ldots,N_{p}-1
\end{equation}

\noindent , $\mathbf{B}^{j} = \left[D_{0}^{j}-\hat{u}_{0}^{j}\left[M_{01}^{j}\right] \quad \ldots \quad  D_{N_{p}-1}^{j}-\hat{u}_{0}^{j}\left[M_{0N_{p}-1}^{j}\right]\right]^{T}$, $\mathbf{X}^{j} = \left[\hat{u}_{0}^{j} \ldots \hat{u}_{N_{p}-1}^{j}\right]^{T}$ and 

\begin{equation}
 \mathbf{A}^{j} = \left[A\right] =  \begin{bmatrix} \left[M_{11}^{j}\right] & \ldots & \left[M_{1N_{p}-1}^{j}\right] \\ \ldots & \ldots & \ldots \\ \rule{0pt}{2.5ex}  \left[M_{N_{p}-1 1}^{j}\right] & \ldots & \rule{0pt}{2.5ex} \left[M_{N_{p}-1N_{p}-1}^{j}\right] \end{bmatrix}
\end{equation}

\noindent Here, the terms $\left[M_{mn}^{j}\right]$ are given by equation \eqref{massMatrix2D}. Then the reconstructed degrees of freedom are given by $\mathbf{X}^{j} = (\mathbf{A}^{j})^{-1}\mathbf{B}^{j}$. This formulation will work for any polynomial basis. Now, we can get the reconstructed polynomial solution in $\mathbf{I}_{j}$ by equation \eqref{modalForm2D}. We call this limiting procedure the compact subcell WENO limiting or CSWENO limiting in short. For solving a system of equations, we use this with a local characteristic field decomposition with the corresponding Jacobians in the $x$ and $y$ directions as explained in \cite{zzsq1}.
\\
\\
\noindent Now, the semi-discrete scheme given in equation \eqref{weakFormScheme2D} along with the limiter is discretized in time by using the TVD Runge-Kutta time discretization introduced in \cite{shu} which is described below briefly.
\\
\\
\noindent If $\{t^{m}\}_{m=0}^{M}$ is a partition of $[0,T]$ and $\Delta t^{m}$ = $t^{m+1}-t^{m}$, $m=0,\ldots,M-1$, the time-marching algorithm reads as follows:
\\
\\
\noindent \textbf{1)} Set $u_{h}^{0}$ = $u_{0h}$, the initial condition;
\\
\\
\textbf{2)} For $m=0,\ldots,M-1$ compute $u_{h}^{m+1}$ from $u_{h}^{m}$ as follows:
\\
$\qquad \qquad$\textbf{(i)} set $u_{h}^{(0)}$ = $u_{h}^{m}$;
\\
$\qquad \qquad$\textbf{(ii)} for $i=1,2,\ldots,d$ compute the intermediate functions:
\begin{equation}\label{timeIntegrate}
u_{h}^{(i)} = \left[\sum_{l=0}^{i-1} \alpha_{il}u_{h}^{(l)} + \beta_{il}\Delta t^{m}L_{h}(u_{h}^{(l)})\right];
\end{equation}

\noindent where $d$ is the order of the time integration.

\noindent \textbf{(iii)} set $u_{h}^{m+1}$ = $u_{h}^{(d)}$;
\\
\\
Some Runge-Kutta time discretization parameters are given in Table \ref{table:tvd}.
\\
\\
\begin{table}
\centering
\begin{tabular}{p{1.5cm}  p{1.2cm}  p{1.2cm}  p{2.5cm} }
\hline
\multicolumn{4}{c}{Runge-Kutta discretization parameters}\\
\hline
order & $\alpha_{il}$ & $\beta_{il}$ & max$\{\alpha_{il}/\beta_{il}\}$ \\
\hline
2 &  1  & 1  & 1 \\
\rule{0pt}{2.5ex}  & $\frac{1}{2}$ $\frac{1}{2}$ & $0$ $\frac{1}{2}$ & \\ [0.5ex]
\hline
3 & 1 & 1 & 1 \\ 
\rule{0pt}{2.5ex}  & $\frac{3}{4}$ $\frac{1}{4}$ &  0 $\frac{1}{4}$ & \\ [0.5ex]
  & $\frac{1}{3}$ 0 $\frac{2}{3}$ & 0 0 $\frac{2}{3}$  &  \\ [0.1ex]
\hline \\ [0.1ex]
\end{tabular}
\caption{TVD Runge-Kutta discretization parameters for orders 2 and 3}
\label{table:tvd}
\end{table}

\noindent We have used the third order TVD Runge-Kutta time discretization for all our calculations.
\\
\\
\section{Results}\label{sec:results}

In this section, we look at some of the results obtained to demonstrate the performance of the limiter (called the compact subcell WENO limiter or CSWENO limiter) described in Section 3. We used Gmsh 4.6.0 software \cite{gmshGR} for the generation of meshes for all our calculations. All the results are obtained using RKDG method and the CSWENO limiter with a third order TVD Runge-Kutta scheme for time integration unless otherwise specified.

\subsection{Accuracy Tests}\label{subsec:accuracyTest}

\begin{figure}[htbp]
\begin{center}
\includegraphics[scale=0.6]{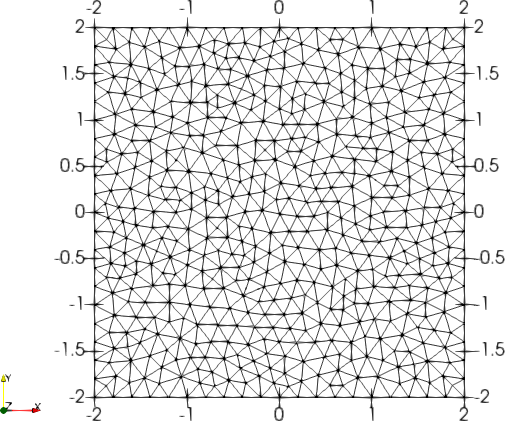}
\caption{Burgers equation - Sample mesh where the mesh points on the boundary are uniformly distributed with cell length $h$ = 4/20}
\label{fig:BurgersSampleMesh}
\end{center}
\end{figure}

\noindent We test the accuracy of the schemes with the CSWENO limiter for scalar and system problems for the two-dimensional test cases. We present the results of the accuracy tests using two-dimensional Burgers equations and the two-dimensional Euler equations. For all the accuracy tests conducted, we have marked all the cells as troubled cells. We used meshes which contain some triangles of irregular size as shown in Figure \ref{fig:BurgersSampleMesh} for all our calculations to illustrate that the limiter retains the order of the scheme even for such meshes.
\\
\\
\noindent \textbf{Example 1:} We solve the two dimensional nonlinear scalar inviscid Burgers equation:
\begin{equation}\label{2dBurgers}
 \frac{\partial u}{\partial t} + \frac{\partial (u^{2}/2)}{\partial x} + \frac{\partial (u^{2}/2)}{\partial y} = 0, \qquad \qquad (x,y) \in [-2,2] \times [-2,2]
\end{equation}
\noindent with the initial condition $u(x,y,0)=0.5+\sin (\pi (x+y)/2)$, with periodic boundary conditions in both directions. The exact solution is smooth till $t=0.5/\pi$. A sample mesh used is shown in Figure \ref{fig:BurgersSampleMesh}. The errors and numerical orders of accuracy are calculated at $t=0.5/\pi$ by marking all the cells as troubled cells and are presented in Table \ref{table:1}. We can see that the CSWENO limiter maintains the order and magnitude of accuracy of the original DG method.
\\
\\
\begin{table}
\centering
\resizebox{\textwidth}{!}{%
\begin{tabular}{|c|c|c|c|c|c|c|c|c|c|}
\hline
\multirow{2}{*}{} &  & \multicolumn{4}{|c|}{DG without limiter} & \multicolumn{4}{|c|}{DG with limiter} \\ \cline{2-10} 
 & $h$ & $L_{1}$ error & Order & $L_{\infty}$ error & Order & $L_{1}$ error & Order & $L_{\infty}$ error & Order \\ \hline
 \multirow{5}{*}{$\mathbf{P}^{1}$} & 4/20 & 2.31E-02 &  & 4.14E-01 &  & 6.64E-03 &  & 8.35E-02 &  \\ \cline{2-10}
  & 4/40 & 4.89E-03 & 2.24 & 8.47E-02 & 2.29 & 1.70E-03 & 1.97 & 2.21E-02 & 1.92 \\ \cline{2-10}
  & 4/80 & 1.01E-03 & 2.27 & 1.84E-02 & 2.20 & 4.37E-04 & 1.96 & 6.30E-03 & 1.81 \\ \cline{2-10}
  & 4/160 & 2.34E-04 & 2.11 & 4.76E-03 & 1.95 & 1.12E-04 & 1.96 & 1.82E-03 & 1.79 \\ \hline
  \multirow{5}{*}{$\mathbf{P}^{2}$} & 4/20 & 4.17E-04 &  & 9.39E-03 &  & 3.97E-04 &  & 9.66E-03 &  \\ \cline{2-10}
  & 4/40 & 5.47E-05 & 2.93 & 1.34E-03 & 2.81 & 5.28E-05 & 2.91 & 1.82E-03 & 2.41 \\ \cline{2-10}
  & 4/80 & 7.59E-06 & 2.85 & 2.08E-04 & 2.69 & 7.32E-06 & 2.85 & 2.84E-04 & 2.68 \\ \cline{2-10}
  & 4/160 & 1.09E-06 & 2.80 & 4.14E-05 & 2.33 & 9.81E-07 & 2.90 & 4.02E-05 & 2.82 \\ \hline
  \multirow{5}{*}{$\mathbf{P}^{3}$} & 4/20 & 3.69E-05 &  & 2.27E-03 &  & 3.87E-05 &  & 8.23E-04 &  \\ \cline{2-10}
  & 4/40 & 2.51E-06 & 3.88 & 1.63E-04 & 3.80 & 2.50E-06 & 3.95 & 6.08E-05 & 3.76 \\ \cline{2-10}
  & 4/80 & 1.73E-07 & 3.86 & 1.13E-05 & 3.85 & 1.63E-07 & 3.94 & 4.16E-06 & 3.87 \\ \cline{2-10}
  & 4/160 & 1.21E-08 & 3.84 & 8.28E-07 & 3.77 & 1.05E-08 & 3.96 & 2.83E-07 & 3.88 \\ \hline
\end{tabular}}
\caption{2D Burgers equation with the initial condition $u(x,y,0)=0.5+\sin (\pi (x+y)/2)$, with periodic boundary conditions in both directions, $t=0.5/\pi$, Triangular mesh with size $h$, $L_{1}$ and $L_{\infty}$ errors for $\mathbf{P}^{1}$, $\mathbf{P}^{2}$ and $\mathbf{P}^{3}$ based DGM}
\label{table:1}
\end{table}

\begin{figure}[htbp]
\begin{center}
\includegraphics[scale=0.6]{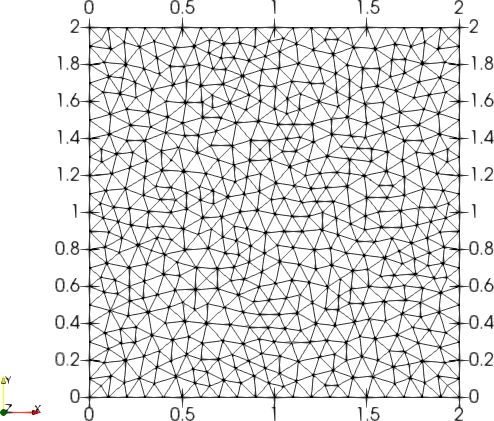}
\caption{2D Euler equations with the initial condition $\rho(x,y,0) = 1+0.2\sin(\pi (x+y))$, $u(x,y,0) = 0.7$, $v(x,y,0) = 0.3$ and $p(x,y,0) = 1.0$ - Sample mesh where the mesh points on the boundary are uniformly distributed with cell length $h$ = 2/20}
\label{fig:EntropyWaveSampleMesh}
\end{center}
\end{figure}

\noindent \textbf{Example 2:} We solve the two dimensional Euler equations:
\begin{equation}\label{2dEulerEquations}
\textbf{U}_{t} + \textbf{f(U)}_{x} + \textbf{g(U)}_{y} = 0 \qquad \qquad (x,y) \in [0,2] \times [0,2]
\end{equation}
\noindent where $\textbf{U} = (\rho, \rho u, \rho v, E)^{T}$, $\textbf{f(U)}=u\textbf{U} + (0, p, 0, pu)^{T}$ and $\textbf{g(U)}=v\textbf{U} + (0, 0, p, pv)^{T}$ with $p = (\gamma -1)(E-\frac{1}{2}\rho (u^{2}+v^{2}))$ and $\gamma = 1.4$. Here, $\rho$ is the density, $(u,v)$ is the velocity, $E$ is the total energy and $p$ is the pressure. The initial conditions are given by $\rho(x,y,0) = 1+0.2\sin(\pi (x+y))$, $u(x,y,0) = 0.7$, $v(x,y,0) = 0.3$ and $p(x,y,0) = 1.0$ and we use periodic boundary conditions in both directions. The exact solution is given by $\rho(x,y,t) = 1+0.2\sin(\pi (x+y-t))$, $u(x,y,t) = 0.7$, $v(x,y,t) = 0.3$ and $p(x,y,t) = 1.0$.
A sample mesh used is shown in Figure \ref{fig:EntropyWaveSampleMesh}. 
The errors in density and numerical orders of accuracy are calculated at $t=2$ by marking all the cells as troubled cells and are presented in Table \ref{table:2}. Again, we can see that the CSWENO limiter maintains the order and magnitude of accuracy of the original DG method.
\\
\\
\begin{table}
\centering
\resizebox{\textwidth}{!}{%
\begin{tabular}{|c|c|c|c|c|c|c|c|c|c|}
\hline
\multirow{2}{*}{} &  & \multicolumn{4}{|c|}{DG without limiter} & \multicolumn{4}{|c|}{DG with limiter} \\ \cline{2-10} 
 & $h$ & $L_{1}$ error & Order & $L_{\infty}$ error & Order & $L_{1}$ error & Order & $L_{\infty}$ error & Order \\ \hline
 \multirow{5}{*}{$\mathbf{P}^{1}$} & 2/20 & 6.69E-03 &  & 2.87E-02 &  & 2.72E-03 &  & 6.95E-03 &  \\ \cline{2-10}
  & 2/40 & 1.51E-03 & 2.15 & 7.58E-03 & 1.92 & 6.90E-04 & 1.98 & 1.81E-03 & 1.94 \\ \cline{2-10}
  & 2/80 & 3.13E-04 & 2.27 & 1.54E-03 & 2.30 & 1.74E-04 & 1.99 & 4.53E-04 & 2.00 \\ \cline{2-10}
  & 2/160 & 8.05E-05 & 1.96 & 3.57E-04 & 2.11 & 4.44E-05 & 1.97 & 1.17E-04 & 1.95 \\ \hline
  \multirow{5}{*}{$\mathbf{P}^{2}$} & 2/20 & 1.02E-04 &  & 2.27E-03 &  & 8.21E-05 &  & 3.31E-03 &  \\ \cline{2-10}
  & 2/40 & 1.15E-05 & 3.15 & 3.11E-04 & 2.87 & 1.03E-05 & 2.99 & 5.97E-04 & 2.47 \\ \cline{2-10}
  & 2/80 & 1.23E-06 & 3.22 & 4.66E-05 & 2.74 & 1.28E-06 & 3.01 & 9.19E-05 & 2.70 \\ \cline{2-10}
  & 2/160 & 1.55E-07 & 2.99 & 6.29E-06 & 2.89 & 1.69E-07 & 2.92 & 1.24E-05 & 2.89 \\ \hline
  \multirow{5}{*}{$\mathbf{P}^{3}$} & 2/20 & 2.31E-06 &  & 4.25E-05 &  & 8.32E-07 &  & 7.24E-06 &  \\ \cline{2-10}
  & 2/40 & 1.50E-07 & 3.95 & 3.25E-06 & 3.71 & 4.92E-08 & 4.08 & 5.31E-07 & 3.77 \\ \cline{2-10}
  & 2/80 & 9.91E-09 & 3.92 & 2.27E-07 & 3.84 & 2.99E-09 & 4.04 & 4.32E-08 & 3.62 \\ \cline{2-10}
  & 2/160 & 7.07E-10 & 3.81 & 1.52E-08 & 3.90 & 1.88E-10 & 3.99 & 2.93E-09 & 3.88 \\ \hline
\end{tabular}}
\caption{2D Euler equations with the initial condition $\rho(x,y,0) = 1+0.2\sin(\pi (x+y))$, $u(x,y,0) = 0.7$, $v(x,y,0) = 0.3$ and $p(x,y,0) = 1.0$, with periodic boundary conditions in both directions, $t=2$, Triangular mesh with size $h$, $L_{1}$ and $L_{\infty}$ errors for density with $\mathbf{P}^{1}$, $\mathbf{P}^{2}$ and $\mathbf{P}^{3}$ based DGM}
\label{table:2}
\end{table}

\begin{figure}[htbp]
\begin{center}
\includegraphics[scale=0.6]{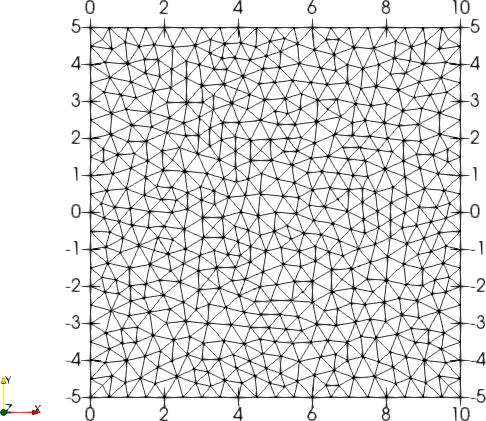}
\caption{2D Euler equations for Isentropic Vortex Problem - Sample mesh where the mesh points on the boundary are uniformly distributed with cell length $h$ = 10/20}
\label{fig:IsenVortexSampleMesh}
\end{center}
\end{figure}

\noindent \textbf{Example 3:} We again solve the two dimensional Euler equations given by \eqref{2dEulerEquations} for the Isentropic Euler Vortex problem suggested by Shu \cite{shu1} as a test case in the domain $[0,10] \times [-5,5]$. The analytical solution is given by: \\ $\rho = \left(1 -  \left(\frac{\gamma - 1}{16\gamma \pi^{2}}\right)\beta^{2} e^{2(1-r^{2})}\right)^{\frac{1}{\gamma-1}}$, $u = 1 - \beta e^{(1-r^{2})} \frac{y-y_{0}}{2\pi}$, $v = \beta e^{(1-r^{2})} \frac{x-x_{0}-t}{2\pi}$, and $p = \rho^{\gamma}$, where $r=\sqrt{(x-x_{0}-t)^{2}+(y-y_{0})^{2}}$, $x_{0}=5$, $y_{0}=0$, $\beta=5$ and $\gamma = 1.4$. We initialize with the analytical solution at $t=0$ and use periodic boundary conditions at the edges of the domain in both directions.
A sample mesh used is shown in Figure \ref{fig:IsenVortexSampleMesh}. 
The errors in density and numerical orders of accuracy are calculated at $t=2$ by marking all the cells as troubled cells and are presented in Table \ref{table:3}. Again, we can see that the CSWENO limiter maintains the order and magnitude of accuracy of the original DG method.
\\
\\

\begin{table}
\centering
\resizebox{\textwidth}{!}{%
\begin{tabular}{|c|c|c|c|c|c|c|c|c|c|}
\hline
\multirow{2}{*}{} &  & \multicolumn{4}{|c|}{DG without limiter} & \multicolumn{4}{|c|}{DG with limiter} \\ \cline{2-10} 
 & $h$ & $L_{1}$ error & Order & $L_{\infty}$ error & Order & $L_{1}$ error & Order & $L_{\infty}$ error & Order \\ \hline
 \multirow{5}{*}{$\mathbf{P}^{1}$} & 10/20 & 4.57E-02 &  & 1.45E-01 &  & 4.35E-02 &  & 1.35E-01 &  \\ \cline{2-10}
  & 10/40 & 1.00E-02 & 2.19 & 3.58E-02 & 2.02 & 1.10E-02 & 1.98 & 3.75E-02 & 1.85 \\ \cline{2-10}
  & 10/80 & 2.27E-03 & 2.14 & 1.11E-02 & 1.69 & 2.81E-03 & 1.97 & 1.01E-02 & 1.90 \\ \cline{2-10}
  & 10/160 & 6.00E-04 & 1.92 & 3.23E-03 & 1.78 & 7.07E-04 & 1.99 & 2.67E-03 & 1.92 \\ \hline
  \multirow{5}{*}{$\mathbf{P}^{2}$} & 10/20 & 6.67E-03 &  & 9.93E-02 &  & 4.89E-03 &  & 8.39E-02 &  \\ \cline{2-10}
  & 10/40 & 7.67E-04 & 3.12 & 1.38E-02 & 2.85 & 6.20E-04 & 2.98 & 1.13E-02 & 2.89 \\ \cline{2-10}
  & 10/80 & 9.26E-05 & 3.05 & 1.88E-03 & 2.88 & 7.91E-05 & 2.97 & 1.80E-03 & 2.65 \\ \cline{2-10}
  & 10/160 & 1.19E-05 & 2.96 & 3.12E-04 & 2.59 & 1.02E-05 & 2.95 & 2.73E-04 & 2.72 \\ \hline
  \multirow{5}{*}{$\mathbf{P}^{3}$} & 10/20 & 8.11E-04 &  & 8.49E-03 &  & 8.27E-04 &  & 7.82E-03 &  \\ \cline{2-10}
  & 10/40 & 5.03E-05 & 4.01 & 6.49E-04 & 3.71 & 5.28E-05 & 3.97 & 5.85E-04 & 3.74 \\ \cline{2-10}
  & 10/80 & 3.35E-06 & 3.91 & 5.28E-05 & 3.62 & 3.44E-06 & 3.94 & 4.53E-05 & 3.69 \\ \cline{2-10}
  & 10/160 & 2.32E-07 & 3.85 & 4.09E-06 & 3.69 & 2.32E-07 & 3.89 & 3.34E-06 & 3.76 \\ \hline
\end{tabular}}
\caption{2D Euler equations for the Isentropic Vortex problem with periodic boundary conditions in both directions, $t=2$, Triangular mesh with size $h$, $L_{1}$ and $L_{\infty}$ errors for density with $\mathbf{P}^{1}$, $\mathbf{P}^{2}$ and $\mathbf{P}^{3}$ based DGM}
\label{table:3}
\end{table}

\subsection{Test Cases With Shocks}\label{subsec:TestCasesShocks}

\begin{figure}[htbp]
  \centering
  \subfloat[Comparison of density solution on $y=0.5$ line of Shock entropy wave Problem for  $\mathbf{P}^{1}$, $\mathbf{P}^{2}$ and $\mathbf{P}^{3}$ based DGM with the exact solution]{\label{fig:ShockDensityWave01}\includegraphics[width=0.45\textwidth]{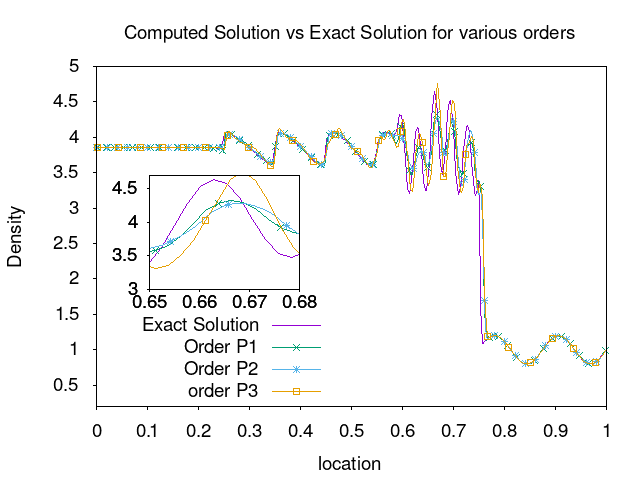}}
  \subfloat[Surface plot for density for $\mathbf{P}^{3}$ based DGM]{\label{fig:ShockDensityWave02}\includegraphics[width=0.45\textwidth]{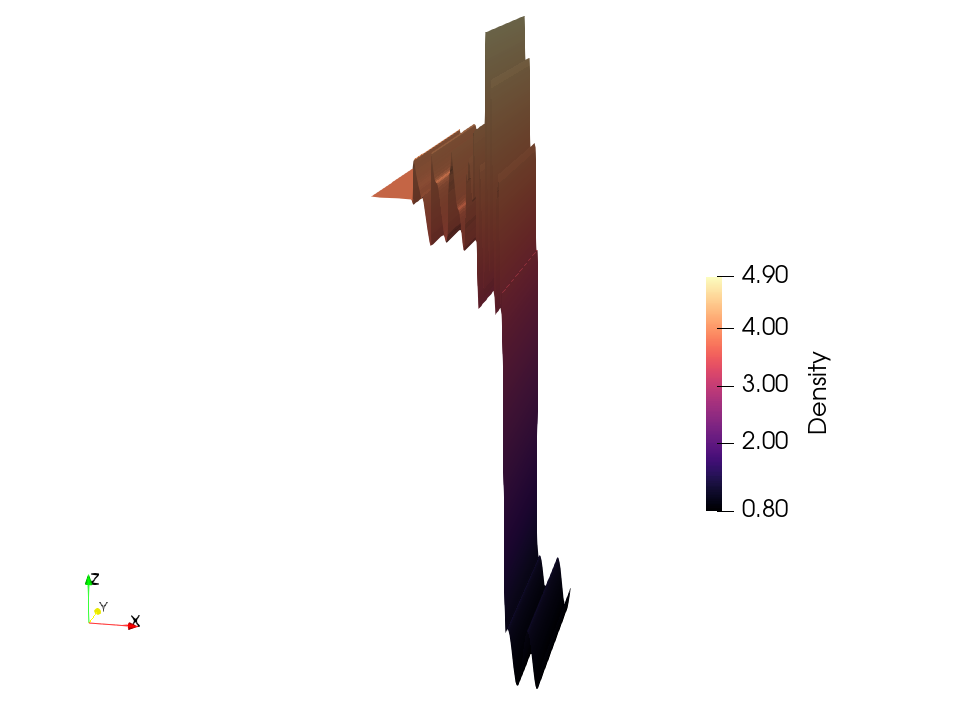}}\hfill
  \caption{Solution of Shock entropy wave Problem at $t=0.178$ with $h=1/200$}
  \label{fig:ShockDensityWave}
\end{figure}

\noindent We now test the CSWENO (compact subcell WENO) limiter for problems with solutions having shocks. 
\\
\\
\noindent \textbf{Example 4:} We solve the problem of shock interaction with entropy waves as proposed in \cite{so2} in the two-dimensional domain. We solve the 2D Euler equations with a moving shock interacting with sine waves in density in the domain $[0,1] \times [0,1]$ with the initial conditions given as $(\rho,u,v,p)=(3.857143,2.629369,0.0,10.333333)$ for $x < 0.125$ and $(\rho,u,v,p)=(1.0+0.2\sin(16\pi x),0,0,1)$ otherwise. Non reflecting boundary condition is applied at $x=0$ and $x=1$ and periodic boundary conditions are applied at the other two boundaries. The computed solution for density obtained at $t=0.178$s using $h=1/200$ at the $y=0.5$ line while using the CSWENO limiter for $\mathbf{P}^{1}$, $\mathbf{P}^{2}$ and $\mathbf{P}^{3}$ based DGM is compared and plotted against the exact solution in Figure \ref{fig:ShockDensityWave01}. A surface color plot for density is shown in Figure \ref{fig:ShockDensityWave02}. We can see that good resolution is obtained in the solution for this problem using CSWENO limiter.
\\
\\

\begin{figure}[htbp]
\begin{center}
\includegraphics[scale=2.0]{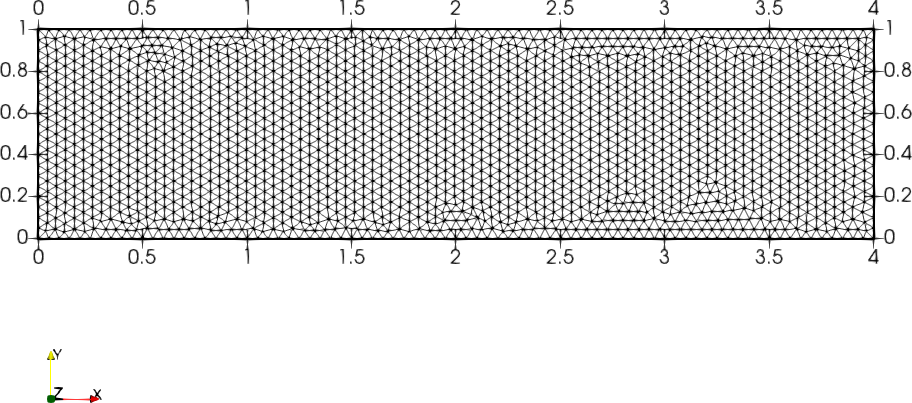}
\caption{Double Mach reflection problem - Sample mesh where the mesh points on the boundary are uniformly distributed with cell length $h$ = 1/20}
\label{fig:DMRSampleMesh}
\end{center}
\end{figure}

\begin{figure}[htbp]
  \centering
  \subfloat[Density contours for the solution at t=0.2 with $\mathbf{P}^{1}$ based DGM]{\label{fig1:P1DGMDMR}\includegraphics[width=0.9\textwidth]{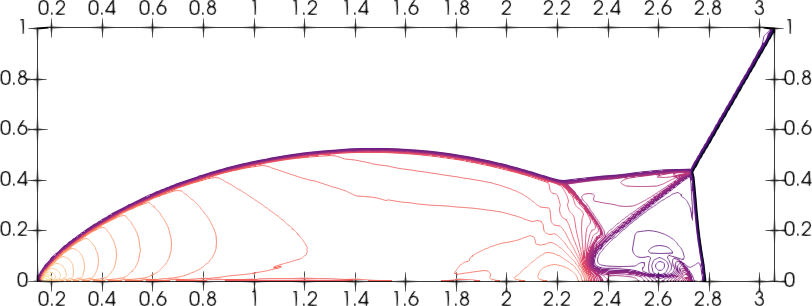}}\hfill
  \subfloat[Density contours for the solution at t=0.2 with $\mathbf{P}^{2}$ based DGM]{\label{fig2:P2DGMDMR}\includegraphics[width=0.9\textwidth]{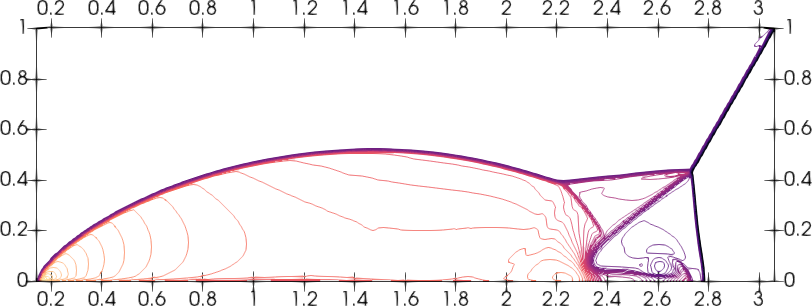}}\hfill
  \subfloat[Density contours for the solution at t=0.2 with $\mathbf{P}^{3}$ based DGM]{\label{fig3:P3DGMDMR}\includegraphics[width=0.9\textwidth]{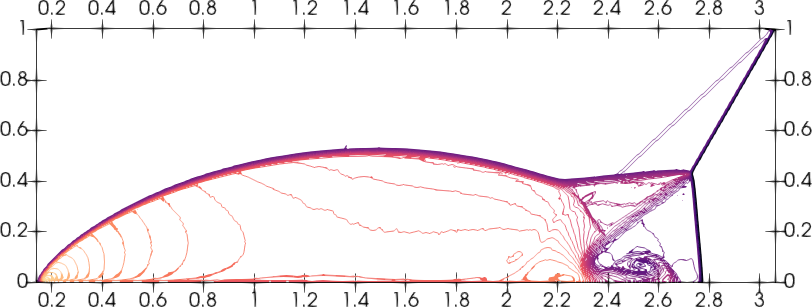}}\hfill
  \subfloat[Density Range]{\label{fig4:DensityRangeDMR}\includegraphics[width=0.9\textwidth]{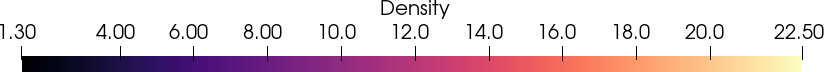}}
  \caption{50 equally spaced density contours for Double Mach reflection solution using the CSWENO limiter with $h=1/200$ for $\mathbf{P}^{1}$(top), $\mathbf{P}^{2}$(middle) and $\mathbf{P}^{3}$(bottom) based DGM}
  \label{fig:CSWENODMRSolution}
\end{figure}

\begin{figure}[htbp]
  \centering
  \subfloat[$\mathbf{P}^{1}$ based DGM]{\label{fig1:P1DGMZoomDMR}\includegraphics[width=0.45\textwidth]{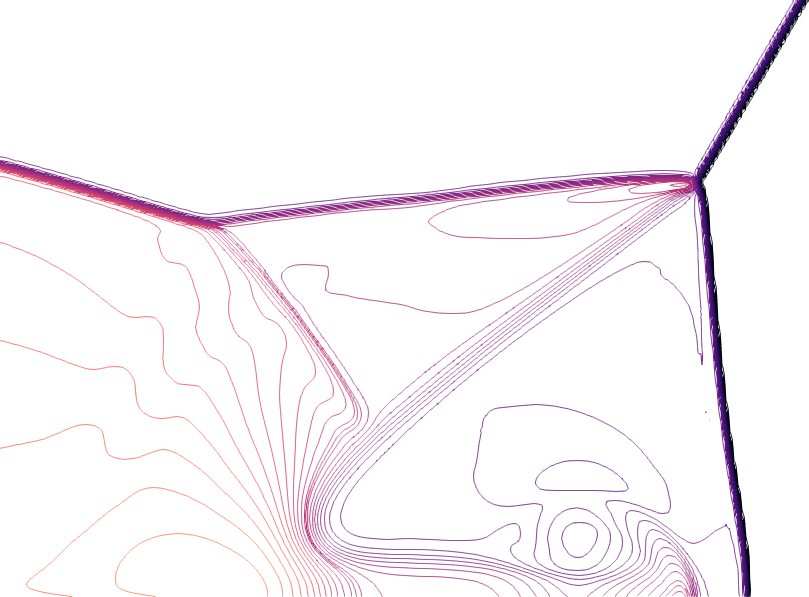}}
  \subfloat[$\mathbf{P}^{2}$ based DGM]{\label{fig2:P2DGMZoomDMR}\includegraphics[width=0.45\textwidth]{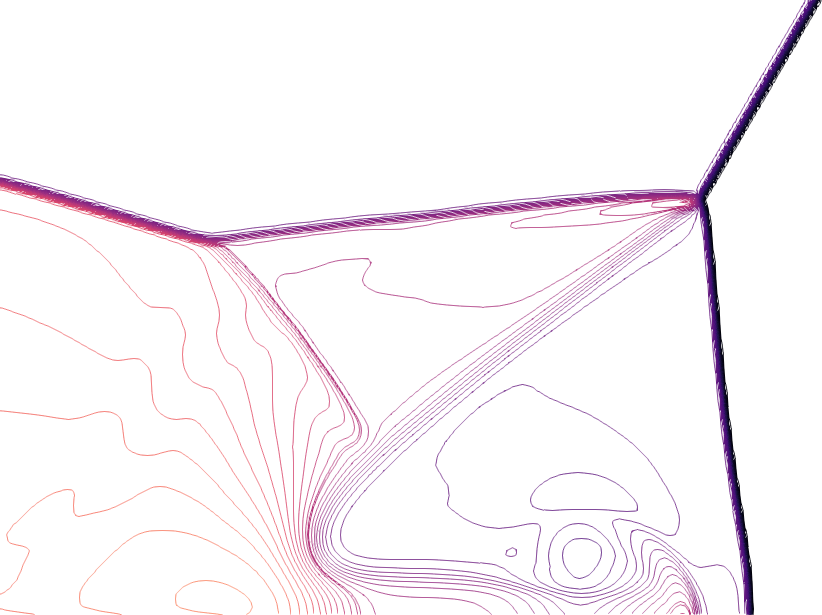}}\hfill
  \subfloat[$\mathbf{P}^{3}$ based DGM]{\label{fig3:P3DGMZoomDMR}\includegraphics[width=0.45\textwidth]{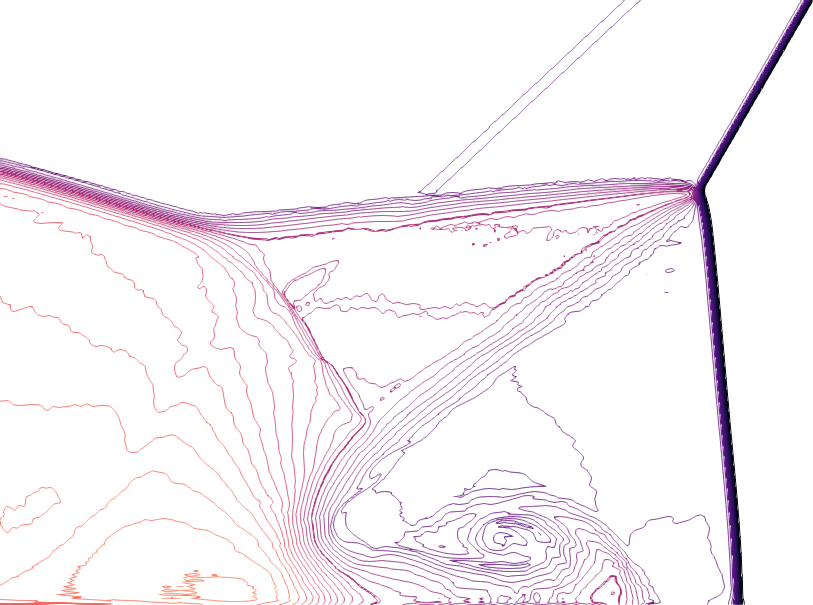}}\hfill
  \subfloat[Density Range]{\label{fig4:DensityRangeCSWENODMR}\includegraphics[width=0.9\textwidth]{DensityRangeDMRUnstructured.png}}
  \caption{Density variation for Double Mach reflection solution using the CSWENO limiter for $\mathbf{P}^{1}$, $\mathbf{P}^{2}$ and $\mathbf{P}^{3}$ based DGM in the region $[2,2.9]\times [0,0.6]$ using 50 equally spaced contours}
  \label{fig:CSWENODMRSolutionZoom}
\end{figure}

\noindent \textbf{Example 5:} We solve the double Mach reflection problem as given in \cite{wc}. We solve the two-dimensional Euler equations in the computational domain $[0,4]\times[0,1]$. Initially, right moving Mach 10 shock is positioned at $x=1/6,y=0$ and it makes an angle $60^{0}$ with the $x$-axis. For the bottom boundary, we impose the exact post shock conditions from $x=0$ to $x=1/6$ and for the rest of the $x$-axis, we use reflective boundary conditions. For the top boundary, we set conditions to describe the exact motion of a Mach 10 shock. We compute the solution upto time $t=0.2$ for a mesh size of $h=1/200$ which contains 370,046 triangles. A sample mesh of size $h=1/20$ is shown in Figure \ref{fig:DMRSampleMesh}. The density contours for the solution obtained using the CSWENO limiter for $\mathbf{P}^{1}$, $\mathbf{P}^{2}$ and $\mathbf{P}^{3}$ based DGM is shown in Figure \ref{fig:CSWENODMRSolution}. A zoom-in view of the density contours near the double Mach stem is shown in Figure \ref{fig:CSWENODMRSolutionZoom} for $\mathbf{P}^{1}$, $\mathbf{P}^{2}$ and $\mathbf{P}^{3}$ based DGM. We can see that the solution obtained using CSWENO limiter is quite well comparable to the solution obtained in \cite{wc}.
\\
\\

\begin{figure}[htbp]
\begin{center}
\includegraphics[scale=2.0]{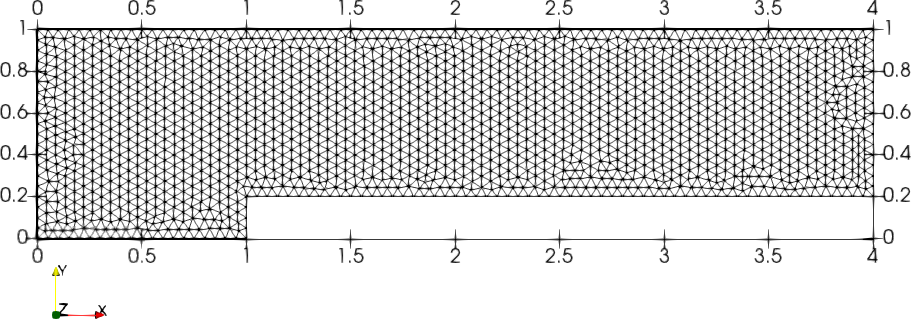}
\caption{Flow over a forward facing step - Sample mesh where the mesh points on the boundary are uniformly distributed with cell length $h$ = 1/20}
\label{fig:ffstepSampleMesh}
\end{center}
\end{figure}

\begin{figure}[htbp]
  \centering
  \subfloat[Density contours for the solution at t=4.0 with $\mathbf{P}^{1}$ based DGM]{\label{fig1:P1DGMffstep}\includegraphics[width=0.9\textwidth]{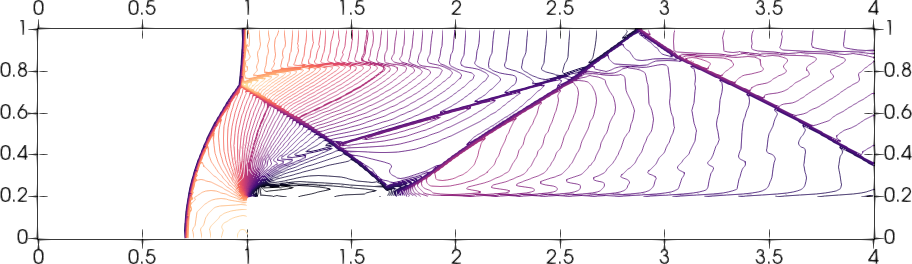}}\hfill
  \subfloat[Density contours for the solution at t=4.0 with $\mathbf{P}^{2}$ based DGM]{\label{fig2:P2DGMffstep}\includegraphics[width=0.9\textwidth]{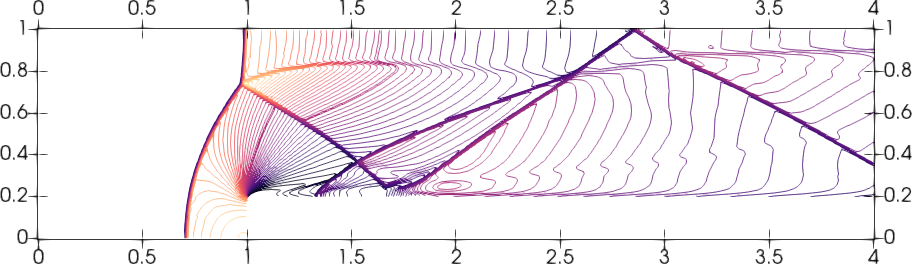}}\hfill
  \subfloat[Density contours for the solution at t=4.0 with $\mathbf{P}^{3}$ based DGM]{\label{fig3:P3DGMffstep}\includegraphics[width=0.9\textwidth]{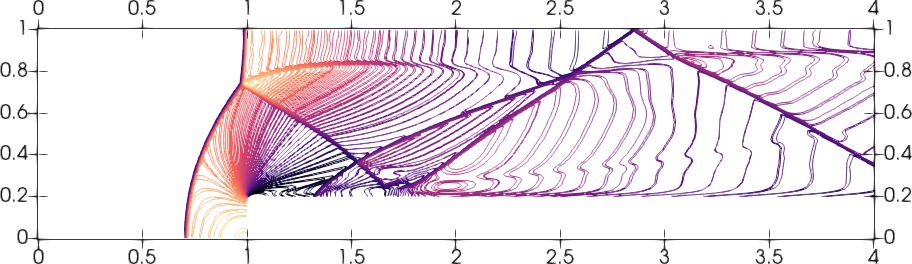}}\hfill
  \subfloat[Density Range]{\label{fig4:DensityRangeffstep}\includegraphics[width=0.9\textwidth]{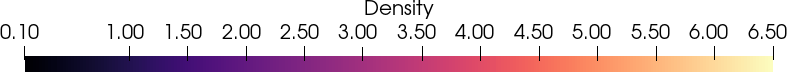}}
  \caption{50 equally spaced density contours for solution of flow over a forward facing step using the CSWENO limiter with $h=1/100$ for $\mathbf{P}^{1}$(top), $\mathbf{P}^{2}$(middle) and $\mathbf{P}^{3}$(bottom) based DGM}
  \label{fig:CSWENOffstepSolution}
\end{figure}

\noindent \textbf{Example 6:} As another test problem, we solve the flow over a forward facing step which is given in \cite{wc}. We solve the two-dimensional Euler equations in a flow set up which contains a right-going Mach 3 uniform flow in a wind tunnel of width 1 unit and length 4 units. The step height is 0.2 units and is located 1 units from the left hand end of the wind tunnel. The problem is initialized by a uniform, right-going Mach 3 flow. Reflective boundary conditions are applied along the walls of the tunnel and in-flow and out-flow boundary conditions are applied at the entrance and the exit, respectively. We compute the solution upto time $t=4.0$ for a mesh size of $h=1/100$ which contains 78982 triangles. A sample mesh of size $h=1/20$ is shown in Figure \ref{fig:ffstepSampleMesh}. The density contours for the solution obtained using the CSWENO limiter for $\mathbf{P}^{1}$, $\mathbf{P}^{2}$ and $\mathbf{P}^{3}$ based DGM is shown in Figure \ref{fig:CSWENOffstepSolution}. We can see that the solution obtained using CSWENO limiter is quite well comparable to the solution obtained in \cite{wc}.
\\
\\

\noindent \textbf{Example 7:} As another test case, we look at the 2D Riemann problem of gas dynamics which is one of the most extensively studied problem which also contains a lot of intricate flow structures. We solve the two-dimensional Euler equations in the domain $[0,1]\times [0,1]$ for 2D Riemann problem configurations (3) and (12) as given by the nomenclature in \cite{ll}. The initial conditions for configurations (3) and (12) are given respectively as

\begin{flalign}\label{2DRiemannInitialConfig03}
 (\rho,u,v,p)(x,y,0) = \begin{cases}
                            (1.5,0,0,1.5) \quad \text{if $x \geq 0.5$ and $y \geq 0.5$} \\
                            (0.5323,1.206,0,0.3) \quad \text{if $x<0.5$ and $y \geq 0.5$} \\
                            (0.138,1.206,1.206,0.029) \quad \text{if $x<0.5$ and $y<0.5$} \\
                            (0.5323,0,1.206,0.3) \quad \text{otherwise}
                        \end{cases}
\end{flalign}

\begin{flalign}\label{2DRiemannInitialConfig12}
 (\rho,u,v,p)(x,y,0) = \begin{cases}
                            (0.5313,0,0,0.4) \quad \text{if $x \geq 0.5$ and $y \geq 0.5$} \\
                            (1,0.7276,0,1) \quad \text{if $x<0.5$ and $y \geq 0.5$} \\
                            (0.8,0,0,1) \quad \text{if $x<0.5$ and $y<0.5$} \\
                            (1,0,0.7276,1) \quad \text{otherwise}
                        \end{cases}
\end{flalign}

\begin{figure}[htbp]
\begin{center}
\includegraphics[scale=0.4]{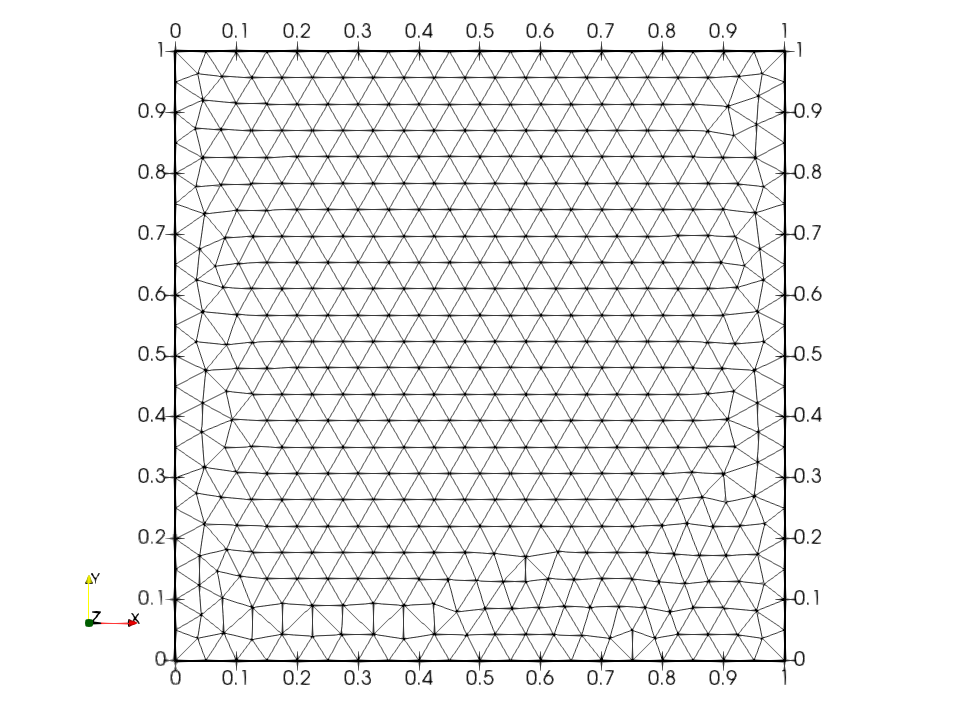}
\caption{2D Riemann problem of gas dynamics - Sample mesh where the mesh points on the boundary are uniformly distributed with cell length $h$ = 1/20}
\label{fig:2dRiemannSampleMesh}
\end{center}
\end{figure}

\begin{figure}[htbp]
  \centering
  \subfloat[Density contours for the solution at t=0.3 with $\mathbf{P}^{1}$ based DGM]{\label{fig1:P1DGM2dRiemannConfig03}\includegraphics[width=0.45\textwidth]{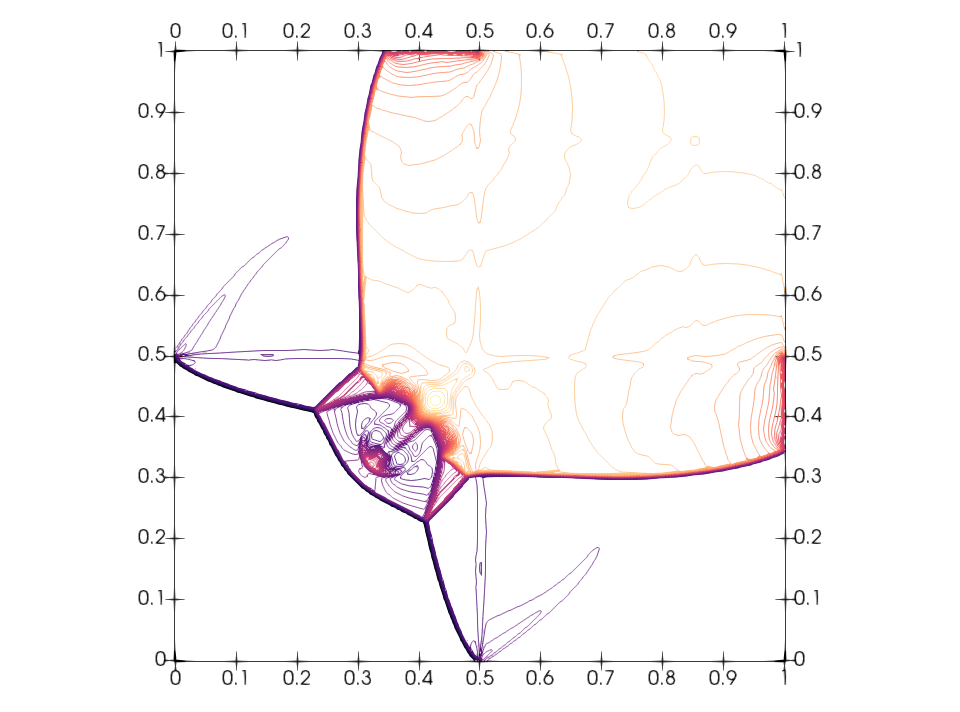}}
  \subfloat[Density contours for the solution at t=0.3 with $\mathbf{P}^{2}$ based DGM]{\label{fig2:P2DGM2dRiemannConfig03}\includegraphics[width=0.45\textwidth]{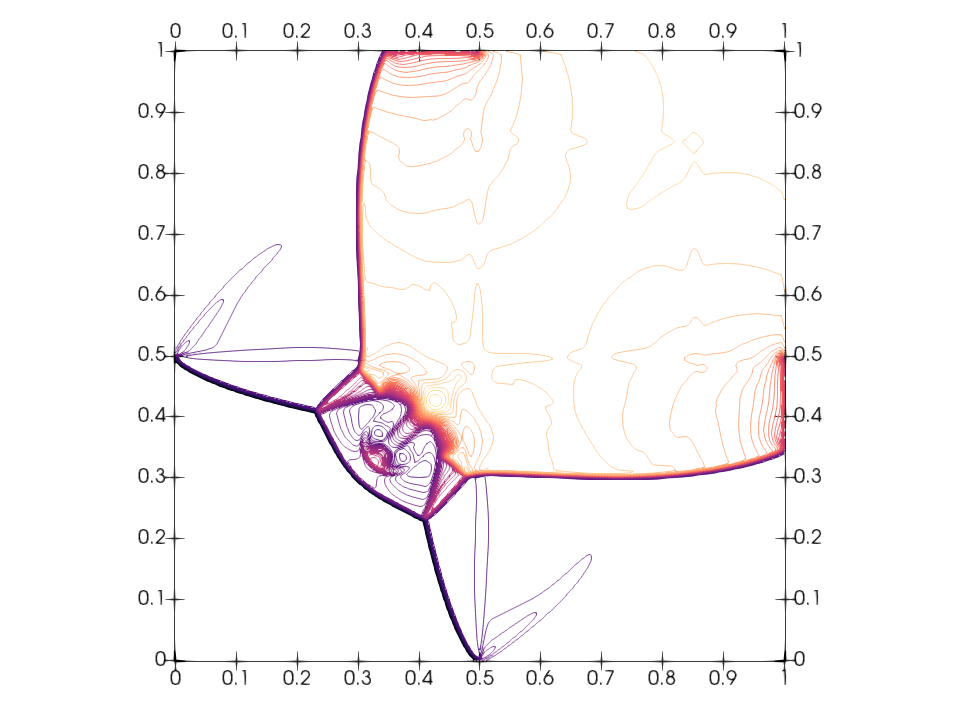}}\hfill
  \subfloat[Density contours for the solution at t=0.3 with $\mathbf{P}^{3}$ based DGM]{\label{fig3:P3DGM2dRiemannConfig03}\includegraphics[width=0.45\textwidth]{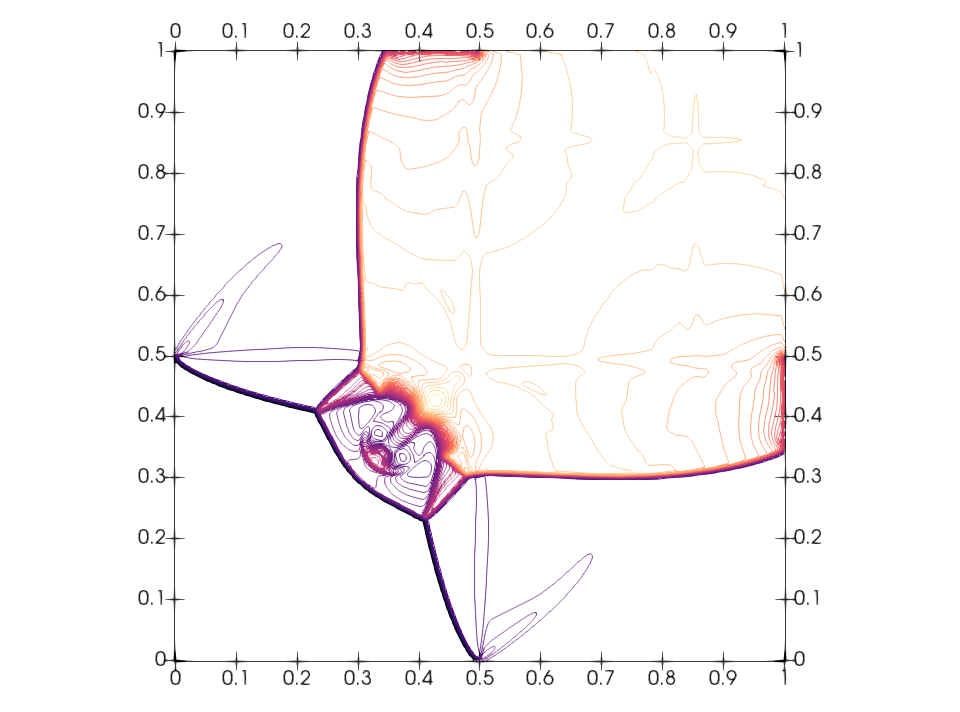}}\hfill
  \subfloat[Density Range]{\label{fig4:DensityRange2dRiemannConfig03}\includegraphics[width=0.9\textwidth]{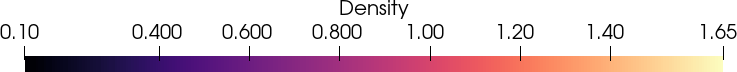}}
  \caption{50 equally spaced density contours for solution at $t=0.3$ for 2D Riemann problem configuration 3 using the CSWENO limiter with $h=1/200$ for $\mathbf{P}^{1}$, $\mathbf{P}^{2}$ and $\mathbf{P}^{3}$ based DGM}
  \label{fig:2dRiemannConfig03Solution}
\end{figure}

\begin{figure}[htbp]
  \centering
  \subfloat[Density contours for the solution at t=0.3 with $\mathbf{P}^{1}$ based DGM]{\label{fig1:P1DGM2dRiemannConfig12}\includegraphics[width=0.45\textwidth]{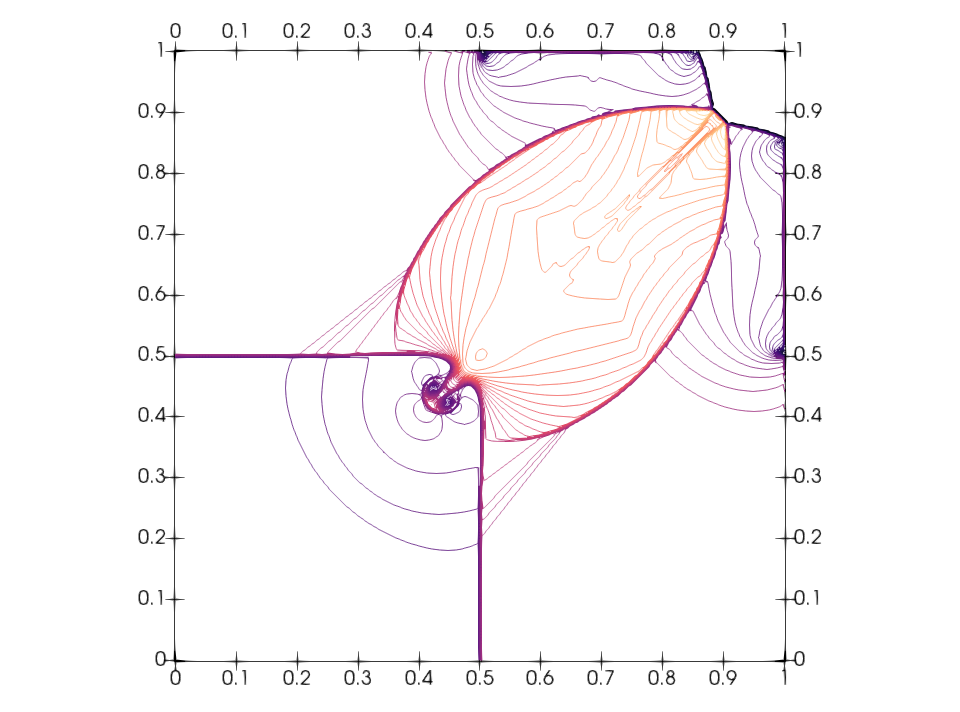}}
  \subfloat[Density contours for the solution at t=0.3 with $\mathbf{P}^{2}$ based DGM]{\label{fig2:P2DGM2dRiemannConfig12}\includegraphics[width=0.45\textwidth]{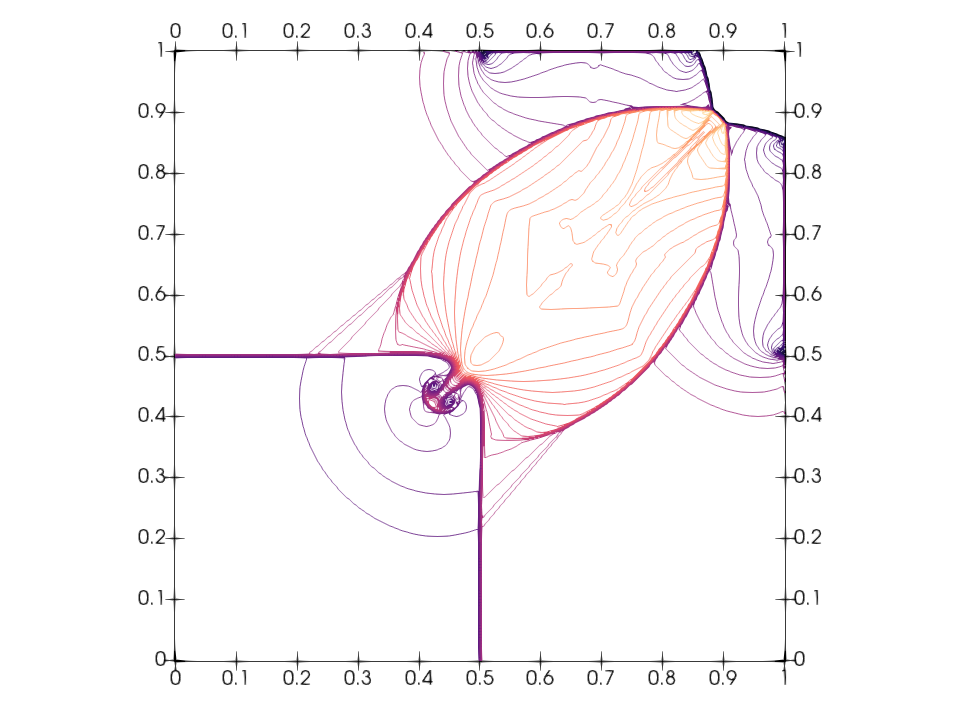}}\hfill
  \subfloat[Density contours for the solution at t=0.3 with $\mathbf{P}^{3}$ based DGM]{\label{fig3:P3DGM2dRiemannConfig12}\includegraphics[width=0.45\textwidth]{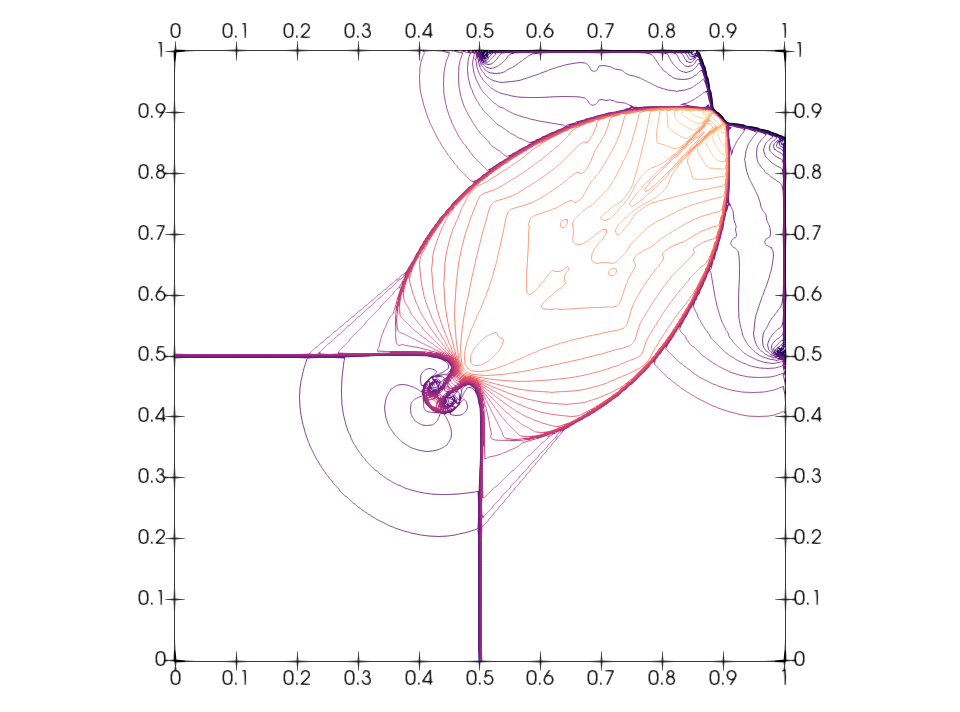}}\hfill
  \subfloat[Density Range]{\label{fig4:DensityRange2dRiemannConfig12}\includegraphics[width=0.9\textwidth]{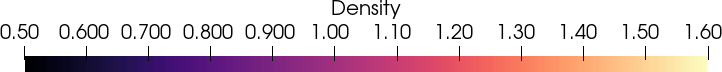}}
  \caption{50 equally spaced density contours for solution at $t=0.25$ for 2D Riemann problem configuration 12 using the CSWENO limiter with $h=1/200$ for $\mathbf{P}^{1}$, $\mathbf{P}^{2}$ and $\mathbf{P}^{3}$ based DGM}
  \label{fig:2dRiemannConfig12Solution}
\end{figure}

We compute the solution upto time $t=0.3$ for configuration (3) and till $t=0.25$ for configuration (12). We use a mesh size of $h=1/200$ which contains 92552 triangles. A sample mesh of size $h=1/20$ is shown in Figure \ref{fig:2dRiemannSampleMesh}. The density contours for the solution obtained using the CSWENO limiter for $\mathbf{P}^{1}$, $\mathbf{P}^{2}$ and $\mathbf{P}^{3}$ based DGM are shown in Figures \ref{fig:2dRiemannConfig03Solution} and \ref{fig:2dRiemannConfig12Solution} respectively for configurations (3) and (12). We can see that the solution obtained using CSWENO limiter is quite well comparable to the solution obtained in \cite{ll}.

\section{Conclusions:}\label{sec:conc}

We have generalized the compact subcell WENO (CSWENO) limiting strategy for the solution of hyperbolic conservation laws using discontinuous Galerkin method proposed in \cite{srspkmr1} to unstructured triangular meshes.~Using this strategy, we identify the troubled cells and use only the immediate neighbors by dividing them into subcells appropriately based on the order of accuracy of the scheme. These new cells are used for the WENO reconstruction.~This formulation is different from the subcell limiting strategy of Dumbser et al \cite{dl1} and Giri et al \cite{gq1} which is much more accurate but quite complicated as they use subcells in an \textit{a posteriori} limiting strategy. This limiting strategy can be used with any of the WENO reconstructions given in \cite{dk1} (called type I WENO reconstruction), or \cite{zqsd} (called type II WENO reconstruction) or \cite{lz1} (mixed reconstruction) or the more recent methods given in \cite{zq1} and \cite{zs4} by dividing the immediate neighbors appropriately. We use the WENO reconstruction given in \cite{zs4} as it is quite simple in implementation and extension to higher orders is easy. We emphasize that our procedure works with any of the WENO reconstructions listed above by dividing the neighbors of the troubled cell as is required by the reconstruction. We termed this limiting procedure as the compact subcell WENO limiter (CSWENO limiter). We have tested the accuracy of this limiter using various standard test cases containing smooth solutions and calculating the numerical order of accuracy. We have also provided numerical results with shocks for standard test cases which are solutions of two-dimensional Euler equations to illustrate the performance of the limiter.


\bibliographystyle{ieeetr}

\bibliography{references}

\end{document}